\documentclass[12pt,leqno]{article}

\usepackage{amsthm}

\usepackage{amsmath}
\usepackage{amscd}
\usepackage{amsmath}
\usepackage{amssymb}
\usepackage{latexsym}
\makeatletter

\@addtoreset{equation}{section}
\makeatother

\newtheorem{thm}{Theorem}[section]
\newtheorem{pr}[thm]{Proposition}
\newtheorem{df}[thm]{Definition}
\newtheorem{lm}[thm]{Lemma}
\newtheorem{cor}[thm]{Corollary}

\newtheorem{rmk}[thm]{Remark}
\newtheorem{ex}[thm]{Example}
\newtheorem{df-lm}[thm]{Definition-Lemma}

\newcommand{\Spec}{{\rm Spec}}

\newcommand{\Ql}{\overline{\mathbb{Q}}_\ell}
\newcommand{\rk}{{\rm rk}}
\newcommand{\sw}{{\rm Sw}}
\newcommand{\dt}{{\rm dimtot}}
\newcommand{\Proj}{\mathbb{P}}

\newcommand{\F}{\mathbb{F}}
\newcommand{\cyc}{{\rm cyc}}

\input xy \xyoption{all} \CompileMatrices

\begin{document}

\title{On continuity of local epsilon factors of $\ell$-adic sheaves}
\author{Daichi Takeuchi
\thanks{Graduated School of Mathematical Science, The University of Tokyo. 
\texttt{Email: daichi@ms.u-tokyo.ac.jp}}}
\date{}
\maketitle

\begin{abstract}
Let $S$ be a noetherian scheme and $f\colon X\to S$ be a smooth morphism of relative dimension $1$. 
For a locally constant sheaf on the complement of a divisor in $X$ flat over $S$, Deligne and Laumon proved that 
the universal local acyclicity is equivalent to the local constancy of Swan conductors. In this article, 
assuming the universal local acyclicity, we show an analogous result of the continuity of local epsilon factors. 
We also give a generalization of this result to a family of isolated singularities. 
\end{abstract}

\section{Introduction}\label{intro}

Let $S$ be a noetherian scheme and $Y$ be a smooth separated curve over $S$. 
In the paper \cite{semi-c}, Deligne-Laumon investigate the relation between 
local acyclicity and ramification along a boundary. 

To be precise, let us fix notations. Let $\Lambda$ be a finite local ring with characteristic 
invertible in $S$. Let $Z$ be a closed subscheme of $Y$ which is finite flat over $S$. 
For a locally constant sheaf $\mathcal{F}$ of finite free $\Lambda$-modules on 
the complement $U=Y\setminus Z$, define a function 
$\varphi_{\rm dt}\colon S\to\mathbb{Z}$ by 
\begin{equation}\label{locfunc}
s\mapsto\sum_{\bar{z}\in Z_{\bar{s}}}{\rm dimtot}_{\bar{z}}\mathcal{F}|_{U_{\bar{s}}}, 
\end{equation}
where $\bar{s}$ is a geometric point above $s$ with algebraically closed residue field and 
${\rm dimtot}_{\bar{z}}=\dim+{\rm Sw}_{\bar{z}}$ is the total dimension function 
with respect to the valuation of the function field of $Y_{\bar{s}}$ defined by the closed point 
$\bar{z}$. The main result in \cite{semi-c} states the following. 
\begin{thm}(\cite[Th\'eor\`eme 2.1.1]{semi-c})\label{scDL}
Let $j\colon U\to Y$ be the inclusion. 
Assume that $S$ is excellent. 
Then the structure morphism 
$Y\to S$ is universally locally acyclic relatively to $j_!\mathcal{F}$
if and only if the function $\varphi_{\rm dt}$ is locally constant, i.e. there exists a unique 
map $\pi_0(S)\to\mathbb{Z}$ which makes the diagram 
\begin{equation*}
\xymatrix{
S\ar[rr]^{\varphi_{\rm dt}}\ar[rd]&&\mathbb{Z}\\&\pi_0(S)\ar[ru]
}
\end{equation*}
commutative. Here $\pi_0(S)$ is the set of connected components and the map 
$S\to\pi_0(S)$ is the canonical one. 
\end{thm}
In particular, the universal local acyclicity implies the local constancy of the function $\varphi_{\rm dt}$. 
In this article, we would like to investigate an arithmetic analogue of this, replacing total dimensions by 
{\it local epsilon factors}. 

Let $k$ be a finite field of characteristic $p$ 
and $T$ be a henselian trait of equal characteristic with residue field $k$. 
Fix a non-trivial additive character $\mathbb{F}_p\to\Ql^\times$, where $\ell$ is a prime number different from $p$. 
For a smooth $\ell$-adic sheaf $V$ on the generic point $\eta$ of $T$, i.e. an $\ell$-adic representation of the 
absolute Galois group of $\eta$, and a non-zero rational $1$-form $\omega\in\Omega^1_\eta$, 
Langlands-Deligne define a constant $\varepsilon_0(T,V,\omega)$ in $\Ql^\times$, the local epsilon factor,  as a candidate which should 
corresponds via local Langlands correspondence to the local constant appearing in the functional equation of
automorphic local L-functions \cite{Del}, \cite{Lau}. In this article, we use this constant to measure ramifications of sheaves. 

In the  papers \cite{Y1} and \cite{Y3}, Yasuda generalizes local epsilon factors to representations in torsion coefficients. 
His theory also can treat the case where the residue field $k$ is a general perfect field of characteristic $p>0$. 
Using his works, our results are valid in this general setting. 
Although local epsilon factors are defined in mixed characteristic cases, we only treat the equal characteristic cases. See the comment 
before the summary of the proof below. 

We give a precise description of our results. Let $S$ be a normal connected scheme of finite type over $\F_p$. 
Let $Y$ be a smooth relative curve over $S$ and $Z$ be a closed subscheme of $Y$ which is finite over $S$. 
Fix a section $\omega\in\Gamma(Y,\Omega^1_{Y/S})$ which generates the relative cotangent sheaf $\Omega^1_{Y/S}$ 
around $Z$. Let $\Lambda$ be a finite local ring with residue characteristic $\ell$ which is invertible in $\F_p$. 
For a locally constant sheaf $\mathcal{F}$ of finite free $\Lambda$-modules on the complement $U=Y\setminus Z$, 
consider the function 
\begin{equation}\label{locepfunc}
\varphi_{\rm ep}\colon |S|\to\Lambda^\times,
\end{equation}
from the set of closed points in $S$, defined by 
\begin{equation*}
s\mapsto\prod_{z\in Z_s}(-1)^{a_z}\varepsilon_{0,\Lambda}(Y_{s(z)},\mathcal{F}|_{U\times_YY_{s(z)}},\omega). 
\end{equation*}
Here $a_z$ is the integer $[k(z)\colon k(s)]{\rm dimtot}_z\mathcal{F}|_{U_s}$ and $\varepsilon_{0,\Lambda}(-)$ 
is the theory of local epsilon factors in torsion coefficients defined by Yasuda \cite{Y1}. 
The symbol $Y_{s(z)}$ means the henselization of $Y_s=Y\times_Ss$ at $z$. 
As an analogue of Theorem \ref{scDL}, we prove that this function satisfies {\it the reciprocity law}.  
\begin{thm}\label{eplocep}
Further assume that the structure map $Y\to S$ is universally locally acyclic relatively to 
the $0$-extension $j_!\mathcal{F}$ for the open immersion $j\colon U\to Y$. 
If, Zariski-locally on $Y$, there exists an \'etale morphism $f\colon Y\to\mathbb{A}^1_S$ such that 
$\omega=df$, there  exists a unique character $\rho_\mathcal{F}\colon 
\pi^{ab}_1(S)\to\Lambda^\times$ which makes the diagram 
\begin{equation*}
\xymatrix{
|S|\ar[rr]^{\varphi_{\rm ep}}\ar[rd]&&\Lambda^\times\\&\pi^{ab}_1(S)\ar[ru]_{\rho_\mathcal{F}}
}
\end{equation*}
commutative. Here the arrow $|S|\to\pi^{ab}_1(S)$ sends closed points $s$ to geometric Frobeniuses ${\rm Frob}_s$. 
\end{thm}
Let $\varphi_{\rm dt}\colon S\to \mathbb{Z}$ be the function defined as in (\ref{locfunc}). By the compatibility of $\varepsilon_{0,\Lambda}$ 
with unramified twists, we have $\rho_{\mathcal{F}(1)}\cdot\rho_{\mathcal{F}}^{-1}=\chi_{\rm cyc}^{\varphi_{\rm dt}}$. 
Here $\chi_{\rm cyc}$ denotes the $\ell$-adic cyclotomic character and $\chi_{\rm cyc}^{\varphi_{\rm dt}}$ is 
the character sending ${\rm Frob}_s\mapsto q_s^{\varphi_{\rm dt}(s)}$, where $q_s$ is the number of the field $k(s)$. 
The symbol $(1)$ means the Tate twist. 

We require that the differential $\omega$ is Zariski-locally of the form $\omega=df$ for an \'etale morphism 
$f\colon Y\to\mathbb{A}^1_S$. This assumption is necessary to apply Theorem \ref{flatan}. 

In \cite{Sai17}, Saito gives a generalization of Theorem \ref{scDL} to a family of isolated singularities. Similarly 
our result also has a generalization to a setting analogous to his as follows. 

Consider the following commutative diagram 
\begin{equation}\label{sas}
\xymatrix{
Z\ar@{^{(}-_>}[r]&X\ar[rd]_g\ar[rr]^{f}&&Y\ar[ld]\ar[r]^t&\mathbb{A}^1_S\ar[lld] \\
&&S 
}
\end{equation}
of $\F_p$-schemes of finite type. We suppose that $S$ is normal of finite type over $\F_p$, that $Z$ is a closed subscheme 
of $X$ finite over $S$, and that $t\colon Y\to\mathbb{A}^1_S$ is \'etale. 
Let $K$ be a complex of constructible sheaves of $\Lambda$-modules with finite tor-dimension on $X$. 
We further assume that $f|_{X\setminus Z}$ and $g$ are universally locally acyclic relatively to $K$. 
For a closed point $s\in S$, the fibers define a sequence of morphisms
\begin{equation}\label{seqcha}
Z_s\to X_s\xrightarrow{f_s}Y_s\xrightarrow{t}\mathbb{A}^1_s. 
\end{equation}
Then the morphism $f_s$  is universally locally acyclic relatively to $K|_{X_s}$ outside $Z_s$, which is finite. 
In this sense, we refer to (\ref{sas}) as a family of isolated singularities relative to $K$. For a point $z\in Z_s$, 
the vanishing cycles complex $R\Phi_{f_s}(K)_z$ is defined from (\ref{seqcha}) and is a complex of 
representations of the absolute Galois group of the generic point of the henselization $Y_{s(z)}$. Thus we can take its 
local epsilon factor $\varepsilon_{0,\Lambda}(Y_{s(z)},R\Phi_{f_s}(K)_z,dt)$. 
We finally consider a function $\theta_{\rm ep}\colon  |S|\to \Lambda^\times$ defined by 
\begin{equation*}
s\mapsto \prod_{z\in Z_s}(-1)^{a_z}\varepsilon_{0,\Lambda}(Y_{s(z)},R\Phi_{f_s}(K)_z,dt) 
\end{equation*}
for $a_z=[k(z)\colon k(s)]{\rm dimtot}_z R\Phi_{f_s}(K)_z$. The following is our main result. 
\begin{thm}(Theorem \ref{flatness})\label{flatan}
Assume that $S$ is normal and connected. Then there exists a unique character 
$\pi^{ab}_1(S)\to\Lambda^\times$ which makes the diagram  
\begin{equation*}
\xymatrix{
|S|\ar[rr]^{\theta_{\rm ep}}\ar[rd]&&\Lambda^\times\\&\pi^{ab}_1(S)\ar[ru]
}
\end{equation*}
commutative. 
\end{thm}
Putting $X=Y$, $f={\rm id}_Y$, and $K=j_!\mathcal{F}$, Theorem \ref{eplocep} is a special case of this theorem. 
To obtain the desired character, first we construct it on the generic points and prove that it extends to a character of the whole 
of $S$ (Lemma \ref{tameobj}). This is the reason why we assume that $S$ is normal. Conversely the assumption of normality is 
sufficient for our result; we prove the theorem for a noetherian normal 
scheme $S$, not necessarily of finite type, over $\F_p$. 

The generalization of Theorem \ref{scDL} by Saito \cite[Proposition 2.16]{Sai17} is a key ingredient for 
constructing characteristic cycles. In \cite{ep_char}, putting coefficients on the irreducible components of 
singular supports, we show that we can define refinements of characteristic cycles, called 
epsilon cycles, and we prove that they give global epsilon factors modulo roots of unity as the intersection numbers 
with the $0$-section in the cotangent bundle. 

It may be natural to ask whether Theorem \ref{flatan} also holds in the mixed characteristic setting. 
One difficulty is to fix a normalization of non-trivial additive characters. Note that fixing an additive  
character of the ring of adeles seems to have nothing to do with our problem. 
Suppose that we could fix it and could obtain the desired character. We cannot expect to get this character 
{\it geometrically}. For example, let $S={\rm Spec}(\mathbb{Z}[1/2])$ and consider the map 
$X=\mathbb{A}^1_S\to Y=\mathbb{A}^1_S$ defined by $x\mapsto x^2$ and the constant sheaf $\Lambda$ on $X$. 
The corresponding character $\rho$ of $G_\mathbb{Q}$ should gives quadratic Gauss sums as the values of Frobeniuses.  
Thus $\rho$ must be pure of weight one, which is not geometric in the sense of Fontaine-Mazur. 

Let us give a summary of the proof. 
In the course of the proof of Theorem \ref{flatan}, we need to treat a family of vanishing cycles complexes. 
To deal with it, we use oriented products of topoi, which give a formalism of vanishing cycles complexes over general 
base schemes. 
In Section \ref{vantop}, we review the notion of them. In particular we give a criterion for vanishing cycles complexes to be locally constant 
(Proposition \ref{ula}). After recalling Laumon's cohomological interpretation of local epsilon factors, which relates 
local epsilon factors and vanishing cycles complexes, 
we prove the main result in Section \ref{mr}. Actually we construct a locally constant constructible complex on some topos 
whose rank gives $\varphi_{\rm dt}$ and whose determinant gives $\theta_{\rm ep}$. See Theorem \ref{flat} for the detail. 

We give notation which we use throughout this paper. 
\begin{itemize}
\item We denote by $G_k$ the absolute Galois group of a field $k$. 
\item We denote by $\chi_\cyc\colon G_k\to\mathbb{Z}_\ell^\times$ 
the $\ell$-adic cyclotomic character. For a finite local ring with residue characteristic $\ell$, 
we write the same letter $\chi_{\rm cyc}$ for the composition $G_k\to\mathbb{Z}_\ell^\times\to\Lambda^\times$. 
\item 
For a finite separable extension $k'/k$ of fields, we denote by ${\rm tr}_{k'/k}\colon G_k^{ab}
\to G_{k'}^{ab}$ the transfer morphism induced  
by the inclusion $G_{k'}\hookrightarrow G_k$. The determinant character of the 
induced representation ${\rm Ind}_{G_{k'}}^{G_k}
1_{G_{k'}}$ of the trivial representation is denoted by $\delta_{k'/k}$. 
\item For a scheme $X$ and its point $x$, $k(x)$ is the 
residue field of $X$ at $x$. 
\item For a finite extension $x'/x$ of the spectra of fields, we denote by 
$\deg(x'/x)$ the degree of the extension. When $x={\rm Spec}(k)$ and 
$x'=\Spec(k')$, we also denote it by $\deg(k'/k)$. 
\item Let $x$ be a geometric point on a scheme $X$. We denote the strict henselization of $X$ at $x$ by $X_{(x)}$. On the other hand, we denote the henselization at a point $x\in X$ by $X_{(x)}$. 
More generally, for a finite separable extension $y$ of $x\in X$, we denote the henselization of $X$ at $y$ by $X_{(y)}$. 
\end{itemize}

\section{Preliminaries on Vanishing topoi (\cite{ori})}\label{vantop}
Let $f\colon X\to S$ and $g\colon Y\to S$ be morphisms of schemes. By abuse of notation,
  we also denote the associated \'etale topoi by $X,Y, S$. 
Otherwise stated, geometric points are assumed to be algebraic geometric points.

For the definition of the oriented product $X\overset\gets\times_{S} Y$, 
 we refer to \cite{ori}. This is the universal object 
  of triples $(T\overset\alpha\to X,T\overset\beta\to Y,g\circ\beta\overset\sigma\to f\circ\alpha)$, where $\alpha$ and $\beta$ are morphisms of topoi, and $\sigma$ is a natural transformation. In particular, if $Y=S$ and $g$ is the identity, there exists a morphism of topoi 
  $\Psi_f\colon X\to X\overset\gets\times_S S$ such that the two triangles in the diagram 
  \begin{equation}
  \label{va}
  \xymatrix{
X\ar[rrdddd]_{{\rm id}}\ar[rrdd]^{\Psi_f}\ar[rrrrdd]^f\\\\&&X\overset\gets\times_S S\ar[rr]^{p_2}\ar[dd]_{p_1}&&S\\\\&&X
}
\end{equation}
are $2$-commutative.

Here are examples of oriented products. 
\begin{ex}
\begin{enumerate}
\item 
Let $X\to S$ be a morphism of schemes and $s$ be a geometric point 
of $S$. Then, the oriented product $s\overset\gets\times_SX$ 
is canonically isomorphic to the \'etale topos of 
$S_{(s)}\times_SX$. Here $S_{(s)}$ is the strict henselization of $S$ at $s$. 
In particular, we have $s\overset\gets\times_SS\cong S_{(s)}$. 
\item
Let $S$ be a Dedekind scheme. Let $s\in S$ be a closed point of codimension 
$1$, i.e. the local ring $\mathcal{O}_{S,s}$ is a discrete valuation ring. Denote by $\eta$ the 
generic point of the henselization $S_{(s)}$. 
For an $s$-scheme $X$, the vanishing topos $X\overset\gets\times_S(S\setminus s)$ is canonically isomorphic to 
the topos $X\times_{S_{(s)}}\eta$, appeared in \cite{SGA7-2}. 
In particular, $s\overset\gets\times_S(S\setminus s)$ is canonically 
isomorphic to the \'etale topos of $\eta$. 
\end{enumerate}
\end{ex}

 A point on the topos $X\overset\gets\times_{S}Y$ can be considered as a triple 
 denoted by $x\gets y$ consisting of a geometric point $x$ of $X$, a geometric point
 $y$ of $Y$ and a specialization $s=f(x)\gets t=g(y)$, i.e. an $S$-morphism $t\to S_{(s)}$.

 Assume that $X$ and $S$ are quasi-compact and quasi-separated. Under this assumption, the topos $X\overset\gets\times_SS$ is coherent, as is explained in \cite[Section 9]{Org}. 
 Let $\Lambda$ be a finite local ring whose residue characteristic is invertible on $S$.
  For a sheaf of sets or a sheaf of $\Lambda$-modules on $X\overset\gets\times_{S} S$,
  we say that it is {\it constructible} if there exist finite partitions $X=\coprod_{\alpha}X_{\alpha}
  , S=\coprod_{\beta}S_{\beta}$ by locally closed constructible subschemes such that 
  its restrictions to $X_{\alpha}\overset\gets\times_{S}S_{\beta}$ are locally constant constructible
  sheaves. This notion of constructibility is the same with the one given 
in \cite[1.9.3]{point}.

  Let $D_{c}^b(-,\Lambda)$ be the full subcategory of $D^{+}(-,\Lambda)$ consisting of bounded complexes
   whose cohomology sheaves are constructible.
   
   Let $f\colon X\to S$ be a morphism of schemes. Then the morphism
   of topoi $\Psi_{f}\colon X\to X\overset\gets\times_{S}S$ defines a functor $R\Psi_f\colon D^{+}(X,\Lambda)\to
   D^{+}(X\overset\gets\times_{S}S,\Lambda)$, which we call the nearby cycles functor. For an object $K$ of $D^{+}(X,\Lambda)$, and a point $x\gets t$ of $X\overset\gets\times_{S}S$, 
   the stalk $R\Psi_{f}K_{(x\gets t)}$ is canonically identified with $R\Gamma(
   X_{(x)}\times_{S_{(f(x))}}S_{(t)}, K)$ (\cite[1.3]{Ill}).
   
   The relation ${\rm id}=p_{1}\circ \Psi_{f}$ defines a natural transformation of functors
    $p_{1}^{\ast}\to R\Psi_{f}$ by adjunction and by the isomorphism ${\rm id}\to p_1\circ\Psi_f$. The cone of this map defines the vanishing cycles functor $R\Phi
    _{f}\colon D^{+}(X,\Lambda)\to D^{+}(X\overset\gets\times_{S}S,\Lambda)$.
    
    We consider a commutative diagram 
    \begin{equation*}
    \xymatrix{
    &X\ar[rr]^{f}\ar[rd]_{p}&&Y\ar[ld]^{g}\\
    &&S
    }
    \end{equation*}
    of schemes. This defines a morphism $\overleftarrow{g}\colon 
     X\overset\gets\times_{Y}Y\to X\overset\gets\times_{S}S$.
     The canonical isomorphism $\overleftarrow{g}\circ \Psi_{f}\to \Psi_{p}$ induces an isomorphism 
$R\overleftarrow{g}_{\ast}\circ R\Psi_{f}\to R\Psi_{p}$. Let $K$ be 
an object of $D^+(X\overset\gets\times_YY)$. We can compute the stalk 
$R\overset\gets g_\ast K_{(x\leftarrow t)}$ of 
$R\overset\gets g_\ast K$
 at a point $x\leftarrow t$ of 
$X\overset\gets\times_SS$ as follows (\cite[Proposition 1.13]{Ill}). 
Let $s$ and $y$ be the images of $x$ in $S$ and $Y$ respectively. 
The topoi $x\overset\gets\times_SS$ and $x
\overset\gets\times_YY$ are canonically isomorphic to $S_{(s)}$ and $
Y_{(y)}$. Under these identifications, 
$R\overset\gets g_\ast K|_{x
\overset\gets\times_SS}$ can be identified by 
$Rg_{(y)_\ast}(K|_{Y_{(y)}})$, where $g_{(y)}$ is the 
map $Y_{(y)}\to S_{(s)}$ induced from $g$. Hence we have 
\begin{equation}\label{stalkpush}
R\overset\gets g_\ast K_{(x\leftarrow t)}\cong
R\Gamma(Y_{(y)}\times_{S_{(s)}}S_{(t)},K).
\end{equation}

A cartesian diagram
     \begin{equation*}
     \xymatrix{
     &X\ar[dd]_{f}
     &&X_{T}\ar[ll]\ar[dd]^{f_{T}}
    \\\\&S&&T\ar[ll]^{i}
     }
     \end{equation*}
    of schemes defines a $2$-commutative diagram
    \begin{equation*}
    \xymatrix{
    &X_{T}\ar[d]^{i}&X_{T}\overset\gets\times_{T}T\ar[l]^{p_1}\ar[d]^{\overleftarrow{i}}
    &X_{T}\ar[l]^{\ \ \ \ \ \Psi_{f_{T}}}\ar[d]^{i}\\
    &X&X\overset\gets\times_{S}S\ar[l]^{p_1}&X\ar[l]^{\ \ \ \ \ \Psi_{f}}}
    \end{equation*}
   and the base change morphisms define a morphism of distinguished triangles 
   \begin{equation*}
\begin{CD}
 \overset{\gets\ast} ip_1^\ast
@>>> \overset{\gets\ast} iR\Psi_f
@>>> \overset{\gets\ast} iR\Phi_f@>>>\\
@V{\simeq}VV @VVV @VVV\\
 p_1^\ast i^\ast
@>>> R\Psi_{f_T}i^\ast
@>>> R\Phi_{f_T}i^\ast@>>>.\\
\end{CD}
   \end{equation*}
   For an object $K$ of $D^{+}(X,\Lambda)$, we say that the formation of $R\Psi_{f}K$ 
   commutes with the base change $T\to S$ if the middle (hence all) vertical arrow is an isomorphism. 
   In particular, if the formation of $R\Psi_f(K)$ commutes with any finite base change $T\to S$, 
   $R\Psi_f(K)_{(x\gets t)}\cong R\Gamma(X_{(x)}\times_{S_{(s)}}S_{(t)},K)$ is canonically isomorphic to 
   $R\Gamma(X_{(x)}\times_{S_{(s)}}t,K)$, taking $T$ as the closure of the image of $t\to S$. 
   Note that, by {\cite[Proposition 2.7.2]{Sai17}}, $f\colon X\to S$ is universally locally acyclic relatively to $K\in D^{+}(X,\Lambda)$ if and 
   only if 
  the map 
$p_1^\ast  i^\ast K\to  R\Psi_{f_T}i^\ast K$ is an isomorphism for any $i\colon T\to S$. Therefore the formation of $R\Psi_{f}K$ commutes with any base change $T\to S$. 

\subsection{Calculation of vanishing cycles complexes}
   \begin{pr}
   \label{pr}(\cite[Proposition 2.8.]{Sai17})
   Let $f\colon X\to S$ be a morphism of finite type of noetherian schemes and $Z\subset X$
   be a closed subscheme which is quasi-finite over $S$. Let $K$ be an object of $
   D_{c}^b(X,\Lambda)$ such that the restriction of $f\colon X\to S$ to the complement $X\setminus Z \to S$
   is universally locally acyclic relatively to the restriction of $K$.
   
   1. $R\Psi_{f}K$ and $R\Phi_{f}K$ are constructible. 
If $K$ is of finite tor-dimension, so are they. Their formations commute with arbitrary base change. $R\Phi_{f}K$
    is supported on $Z\overset\gets\times_{S}S$.
    
    2. Let $x$ be a geometric point of $X$ and $s=f(x)$ be the geometric point of 
    $S$ defined by the image of $x$ by $f$. Let $t$ and $u$ be geometric points
     of $S_{(s)}$ and $t\gets u$ be a specialization. Then, there exists a distinguished triangle
     \begin{equation*}
     R\Psi_{f}K_{(x\gets t)}\longrightarrow
     R\Psi_{f}K_{(x\gets u)}\longrightarrow
     \bigoplus_{z\in Z_{(x)}\times_{S_{(s)}t}}
     R\Phi_{f}K_{(z\gets u)}\longrightarrow, 
     \end{equation*}
where the first map is the cospecialization and the second one is 
the direct sum of the compositions 
of the maps $R\Psi_{f}K_{(x\gets u)}\to R\Phi_{f}K_{(x\gets u)}$ and the cospecializations 
$R\Phi_{f}K_{(x\gets u)}\to 
R\Phi_{f}K_{(z\gets u)}$.
     \end{pr}
     \proof{
     1.
     The commutativity with base change is proved in {\cite[Proposition 6.1.]{Org}}, taking a compactification and using the proper base change theorem. The last assertion follows from 
     the remark made before this proposition. The constructibility follows from {\cite[Th\'eor\`eme 8.1.]{Org}} and the commutativity with base change. The finiteness of tor-dimension follows from the finiteness of 
cohomological dimension of $R\Psi_f$ {\cite[Proposition 3.1.]{Org}}.

     2. We may assume that $S=S_{(s)}$, that $Z$ is local and finite over $S$, and
      that $X$ is separated over $S$. Let us denote the morphisms obtained by the base change $X\to S$ by the same letter, by abuse of notation.
  By 1, we may further replace $S$ by the normalization in $u$. 
      Consider the following diagram
      \begin{equation*}
      \xymatrix{
     &s\ar[dd]_{i_s}&& t\ar[dd]^{i_t}&&
     \\\\
     &S&&S_{(t)}\ar[ll]_{j}&& u\ar[ll]_{k}
     }
     \end{equation*}
     and define a complex $\Delta$ on $X\times_{S}S_{(t)}$ fitting in the distinguished triangle
     $j^{\ast}K\to Rk_{\ast}(jk)^{\ast}K\to \Delta\to$. 
Note that, since $u$ is the generic point of $S$, the geometric point $t$ can 
be identified with its image in $S$ and that $S_{(t)}$ can be identified 
with the localization of $S$ at $t$. Hence 
the complex $\Delta$ is supported on $Z\times_S
S_{(t)}$. Indeed, let $t'$ be a geometric point of $S_{(t)}$. 
The restriction of $\Delta$ to $X\times_SS_{(t')}$ equals to 
the complex defined in the same way as $\Delta$ replacing $t$ by $t'$. 
Thus the restriction $\Delta|_{X_{t'}}$ is supported on $Z_{t'}$. 
Since the formation of 
     nearby cycles functor commutes with any base change by 1, $(Rj_{\ast}j^{\ast}K)_{x}$ and 
     $(R(jk)_{\ast}(jk)^{\ast}K)_{x}$ are isomorphic to
     $R\Psi_{f}K_{(x\gets t)}$ and  $R\Psi_{f}K_{(x\gets u)}$
      respectively. 
We have isomorphisms
\begin{equation*}
(Rj_\ast\Delta)_{x}\to R\Gamma(Z,Rj_\ast\Delta)\to
R\Gamma(Z\times_SS_{(t)},\Delta). 
\end{equation*}
Since the last term is isomorphic to $\bigoplus_{z\in Z\times_{S}t}
 R\Phi_{f}K_{(z\gets u)}$, the assertion follows. 
      \qed}
      \vspace{0.1in}
      
      To deduce corollaries, we need a following lemma. 
      Let $f\colon X\to S$ be a morphism of schemes. 
     For a (usual) point $x\in X$, we denote by $X_{x}$ the 
     spectrum of the localization $\mathcal{O}_{X,x}$. 
      For a complex $K\in D^+(X,\Lambda)$, 
      we define ${\rm ULA}(K,f)\subset X$ to be the subset consisting of points 
      $x\in X$ such that the morphism $X_{x}\to S$ 
      is universally locally acyclic relatively to $K|_{X_x}$. Note that, 
      for a morphism $g\colon S'\to S$ of schemes, 
      the inverse image of ${\rm ULA}(K,f)$ by $g_{X}\colon X_{S'}:=X\times_SS'
      \to X$ is contained in ${\rm ULA}(g_{X}^\ast K,f_{S'})$, where $f_{S'}$ 
      is the base change 
      $X_{S'}\to S'$. 
      \begin{lm}\label{locconstr}
      Let the notation be as above. 
      \begin{enumerate}
      \item Let $x\to X$ be a geometric point and let $j\colon X_{(x)}\to X$  be 
      the morphism from the strict henselization. Then, we 
      have $j^{-1}({\rm ULA}(K,f))={\rm ULA}(K|_{X_{(x)}},fj)$. 
      \item Assume that 
      $f\colon X\to S$ is of finite type and that $S$ is (hence also $X$ is) 
      noetherian. When $K$ is constructible, 
      ${\rm ULA}(K,f)$ is an open subset of $X$. 
      \end{enumerate}
      \end{lm}
      \proof{
      
      1. 
      Let $y'\in X_{(x)}$ be a point and $y=j(y')\in X$ be its image. 
      Fix a geometric point $\bar{y}\to X_{(x)}$ over $y'$. We also regard $\bar{y}$ as a 
      geometric point of $X$ by $j$. 
      The strict henselizations 
      $X_{(x)(\bar{y})}$ and $X_{(\bar{y})}$ are canonically isomorphic. Since the 
      universal local acyclicity of $X_{(x)y'}:={\rm Spec}(\mathcal{O}_{X_{(x)},y'})\to S$ (resp. $X_y\to S$) is equivalent to that of the 
      morphism $X_{(x)(\bar{y})}\to S$ (resp. $X_{(\bar{y})}\to S$), the assertion follows. 
      
      2. 
      First we show the assertion assuming that the formation of 
      $R\Phi_f(K)$ commutes with arbitrary 
      base change $S'\to S$. Under this assumption, $R\Phi_f(K)$ is 
      constructible by \cite[Th\'eor\`eme 8.1]{Org}. 
     
      Let $x\in X$ be a point in ${\rm ULA}(K,f)$. 
      We need to find an open neighborhood $U$ of $x$ contained 
      in ${\rm ULA}(K,f)$. Since $R\Phi_f(K)$ is a constructible complex on 
      $X\overset\gets\times_SS$ and the restriction 
      of $R\Phi_f(K)$ to $X_x\overset\gets\times_SS$ is acyclic, 
      there is an open neighborhood $U$ of $x$ such that 
      $R\Phi_f(K)$ is acyclic on $U\overset\gets\times_SS$. 
      Hence $f|_U\colon U\to S$ is universally locally acyclic 
     by the commutativity of the formation of $R\Phi_f(K)$. See \cite[Proposition 2.7.2]{Sai17} for the detail. 
      
      Next we consider the general case. 
      By \cite[Th\'eor\`emes 2.1, 8.1]{Org}, there exists a proper surjective morphism  
      $g\colon T\to S$ such that the vanishing cycles complex $R\Phi_{f_T}(K|_{X_T})$ of $K|_{X_T}$ with respect to the base change 
$f_T\colon X_T:=X\times_ST\to T$ is constructible and its formation commutes 
with arbitrary base change $T'\to T$. Let $g_X\colon X_T\to X$ be 
the projection. We already know that the complement 
$Z:=X_T\setminus {\rm ULA}(g_X^\ast K,f_T)$ is closed. Let $Z':=g_X(Z)$ 
be the image, which is closed. 
We show the equality ${\rm ULA}(K,f)=X\setminus Z'=:U'$. 
The intersection ${\rm ULA}(K,f)
\cap Z'$ is empty since $g_X^{-1}({\rm ULA}(K,f))\subset 
{\rm ULA}(g_X^\ast K,f_T)$. Hence it suffices to show that $U'\to S$ 
is universally locally acyclic relatively to $K$. 
Since the base change $U'_T$ is contained in ${\rm ULA}(g_X^\ast K,f_T)$, 
the assertion follows from the oriented cohomological descent \cite[Lemma 6.1.]{LZ}. 
\qed
}

\begin{cor}\label{prcor}
Let $f\colon X\to S$ be a morphism of finite type of noetherian schemes. 
Let $Z\subset X$ be a closed subscheme quasi-finite over $S$. 
Let $x\to X$ be a geometric point and let $s:=f(x)$ be the 
geometric point of $S$ induced from $x$. 

Let $K\in D^b_c(X_{(x)},\Lambda)$ be a constructible complex. 
If $f_{(x)}\colon X_{(x)}\to S_{(s)}$ is universally locally acyclic relatively to $K$ outside $Z\times_XX_{(x)}$, 
$R\Psi_{f_{(x)}}(K)$ is constructible and its formation commutes with arbitrary 
base change $S'\to S_{(s)}$. 
\end{cor}
\proof{
Replacing $X\to S$ by $X\times_SS_{(s)}\to S_{(s)}$, 
we may assume that $S=S_{(s)}$. 
Since $K$ is constructible, after replacing $X$ by an \'etale neighborhood 
of $x$, we may assume that there exists a constructible complex 
$K'$ on $X$ which restricts to $K$. Define $U:={\rm ULA}(K',f)$. 
This is an open subset of $X$ by Lemma \ref{locconstr}.2. By the assumption and Lemma \ref{locconstr}.1, 
we have $(Z\times_XX_{(x)})\cup (U\times_XX_{(x)})=X_{(x)}$. Hence, further replacing $X$, 
we may assume that $Z\cup U=X$. Then, the assertion follows 
from Proposition \ref{pr}.1. 
\qed
}
\vspace{0.1in}

We give a partial generalization of Corollary \ref{prcor} to the following setting. 

Consider a commutative diagram 
\begin{equation}\label{commsc}
\xymatrix{
Y\ar[r]^{g}\ar[d]&T\ar[d]\\
X\ar[r]^f&S
}
\end{equation}
of schemes, to which we refer as a pair $\bold{f}=(f,g)$. 
This induces a morphism of topoi 
$\overset\gets{\bold{f}}\colon Y\overset\gets\times_{X}X\to 
T\overset\gets\times_{S}S$. For a geometric point 
$y$ of $Y$, let $t$ denote the image $g(y)$. This is a 
geometric point of $T$. They are also regarded as geometric points 
of $X$ and $S$ via the vertical arrows of (\ref{commsc}). 
Thus, they give a morphism of topoi 
$f_{(y)}\colon X_{(y)}\to S_{(t)}$, which is also obtained 
from $\overset\gets{\bold{f}}$ via the identifications $X_{(y)}
\cong 
y\overset\gets\times_{X}X$ and $S_{(t)}\cong 
t\overset\gets\times_{S}S$. 

Let $\bold{i}:=(i_T\colon T'\to T,i_S\colon S'\to S)$ be a pair of morphisms of schemes such that the diagram 
\begin{equation}\label{commTS}
\xymatrix{
T'\ar[r]^{i_T}\ar[d]&T\ar[d]\\S'\ar[r]^{i_S}&S
}
\end{equation}
is commutative. The pull-backs give a commutative diagram 
\begin{equation*}
\xymatrix{Y':=Y\times_TT'\ar[r]^{\ \ \ \ \ \ g'}\ar[d]&T'\ar[d]\\X':=X\times_SS'\ar[r]^{\ \ \ \ \ \ f'}&S'
}
\end{equation*}
of schemes, which has a projection to the diagram (\ref{commsc}). Hence we get a commutative diagram 
\begin{equation}\label{commori}
\xymatrix{
Y'\overset\gets\times_{X'}X'\ar[r]^{\overset\gets{\bold{f}'}}\ar[d]^{\overset\gets{\bold{i}}}&T'
\overset\gets\times_{S'}S'\ar[d]^{\overset\gets{\bold{i}}}
\\Y\overset\gets\times_{X}X\ar[r]^{\overset\gets{\bold{f}}}&T\overset\gets\times_{S}S
}
\end{equation}
of topoi. 
\begin{df}\label{relbcdef}
Let $K\in D^+(Y\overset\gets\times_{X}
X,\Lambda)$ be a bounded below complex. 
\begin{enumerate}
\item Let $Z\subset X$ be a closed subset with complement $U$. 
 We say that the pair $\bold{f}=(f,g)$ is (resp. universally) locally acyclic relatively 
to $K$ outside $Z$ if, for every geometric point $y$ of $Y$ with the image $t=g(y)$, 
the morphism $f_{(y)}|_{X_{(y)}\times_XU}\colon X_{(y)}\times_XU\to S_{(t)}$ is (resp. universally) 
locally acyclic relatively to $K|_{X_{(y)}\times_XU}$. 
\item Let $\bold{i}:=(i_T\colon T'\to T,i_S\colon S'\to S)$ be as in (\ref{commTS}). 
We say that the formation of $R\overset\gets{\bold{f}}_\ast K$ 
commutes with a base change $\bold{i}$ if 
the base change map $ \overset\gets{\bold{i}}^\ast R\overset\gets{\bold{f}}_\ast K\to 
R\overset\gets{\bold{f}'}_{\ast} \overset\gets{\bold{i}}^\ast K$ defined from (\ref{commori}) is an isomorphism. 
\end{enumerate}
\end{df}

\begin{lm}\label{finfact}
Assume that $g$ is finite and that $f$ is the identity $S\to S$, and consider the morphism 
of topoi $\overset\gets{\bold{f}}\colon Y\overset\gets\times_SS\to T\overset\gets\times_SS$. 
\begin{enumerate}
\item The formation of $\overset\gets{\bold{f}}_\ast$ commutes with 
base change. 
\item The cohomological dimension of $\overset\gets{\bold{f}}_\ast$ is zero. 
\item Assume that $Y,T,S$ are quasi-separated and quasi-compact. Then the pushforward 
$\overset\gets{\bold{f}}_\ast\colon D^+(Y\overset\gets\times_SS,\Lambda)\to D^+(T\overset\gets\times_SS,\Lambda)$ 
preserves the constructibility. 
\end{enumerate}
\end{lm}
\proof{
1, 2. 
 Let $K$ be a bounded below complex on $Y\overset\gets\times_SS$. Since 
$g$ is finite, for a geometric point $t$ of $T$, 
the restriction of $R\overset\gets{\bold{f}}_\ast K$ to $t\overset\gets\times_SS$ 
is canonically isomorphic to $Rg_{(t)\ast}K$, where 
$g_{(t)}\colon\coprod_{y\in Y\times_{T}t}S_{(y)}\to S_{(t)}$. 
Hence the assertions 1, 2 follow. 

3. Taking a finite stratification of $T$, we may assume that 
$g$ is finite \'etale by 1. Further replacing $T$ by a finite \'etale cover, 
we may assume that $Y$ is isomorphic to the disjoint union of finite copies of $T$. Then the assertion is clear. 
\qed}

\begin{cor}\label{relbc}
Consider a commutative diagram (\ref{commsc}) of schemes. 
Let $Z\subset X$ be a closed subset which is quasi-finite over $S
$.  
Assume that $S$ is noetherian and that $f$ is of finite type. 
Further assume that $g$ is finite. 

Let $K\in D^b_c(Y\overset\gets\times_{X}X,\Lambda)$ 
be a constructible complex. 
Assume that $\bold{f}=(f,g)$ is universally locally acyclic relatively to $K$ outside 
$Z$. 
Then, 
the formation of $R\overset\gets{\bold{f}}_\ast K$ commutes 
with base change $(T'\to T,S'\to S)$. 
\end{cor}
\proof{ 
The morphism $\overset\gets{\bold{f}}$ decomposes into 
\begin{equation*}
Y\overset\gets\times_{X}X\xrightarrow{\overset\gets{f}}Y\overset\gets\times_{S}S\xrightarrow{\overset\gets{g}} T\overset\gets\times_{S}S. 
\end{equation*}
 
The formation of $\overset\gets{g}_\ast$ commutes with base change by Lemma \ref{finfact}. 

Let $\bold{i}=(i_Y\colon Y'\to Y,i_S\colon S'\to S)$ be a pair of morphisms of schemes such that the diagram 
\begin{equation*}
\xymatrix{
Y'\ar[r]^{i_Y}\ar[d]&Y\ar[d]\\S'\ar[r]^{i_S}&S
}
\end{equation*}
is commutative. Let $(y'\gets u')$ be a geometric point of $Y'\overset\gets\times_{S'}S'$ and $(y\gets u)$ be its image in 
$Y\overset\gets\times_SS$. Let $t$ (resp. $t'$) be the geometric point of $S$ (resp. $S'$) defined by 
$y\to Y\to S$ (resp. $y'\to Y'\to S'$). By (\ref{stalkpush}), we have 
\begin{equation*}
(R\overset\gets{f}_\ast K)_{(y\gets u)}\cong R\Gamma(X_{(y)}\times_{S_{(t)}}S_{(u)}, 
K)\cong R\Psi_{f_{(y)}}K_{(y\gets u)}. 
\end{equation*}
Here $f_{(y)}$ is the morphism $X_{(y)}\to S_{(t)}$ of the strict henselizations. 
By Corollary \ref{prcor}, the pull-back of $R\Psi_{f_{(y)}}K$ to $(X_{(y)}\times_{S_{(t)}}S'_{(t')})\overset\gets\times_{S'_{(t')}}S'_{(t')}$ 
is isomorphic to $R\Psi_{f'_{(y)}}(K|_{X_{(y)}\times_{S_{(t)}}S'_{(t')}})$, where $f'_{(y)}\colon X_{(y)}\times_{S_{(t)}}S'_{(t')}\to S'_{(t')}$ is the 
base change of $f_{(y)}$. The assertion follows. 
\qed
}

\medskip
 \begin{pr}\label{ula}
 Let $X,Z,S$ be as in Proposition \ref{pr}. Further assume that $Z\to S$ is finite. 
 Let $Z'$ be the image $f(Z)$. 
 Let $g\colon Z\to Z'$ be the restriction of $f$ and 
$\overset\gets g\colon Z\overset\gets\times_S(S\setminus Z')\to Z'
\overset\gets\times_S(S\setminus Z')$ be the induced morphism of topoi. 
  Let $K$ be an object of $D_{c}^b(X,\Lambda)$ such that 
 the restriction of $K$ to $f^{-1}(Z')$ is acyclic. Assume that 
 $f\colon X\to S$ is universally locally acyclic relatively to $K$ outside $Z$. 
 Then, the following hold. 
\begin{enumerate}
\item The formation 
of $\overset\gets g_\ast(R\Psi_f(K)
|_{Z\overset\gets\times_S(S\setminus Z')})$ commutes with arbitrary base change $S'\to S$. 
\item
The complex $\overset\gets g_\ast(R\Psi_f(K)
|_{Z\overset\gets\times_S(S\setminus Z')})$ is locally constant constructible. 
If $K$ is of finite tor-dimension, so is $\overset\gets g_\ast(R\Psi_f(K)
|_{Z\overset\gets\times_S(S\setminus Z')})$. 
\end{enumerate}
 \end{pr}
 \proof{
 1. 
The formation of $R\Psi_f(K)$ commutes with base change by Proposition \ref{pr}.1. The formation of $\overset\gets{g}_\ast$ commutes with 
base change by Lemma \ref{finfact}.1. 

2. 
The constructibility of $R\Psi_f(K)
|_{Z\overset\gets\times_S(S\setminus Z')}$ follows from Proposition \ref{pr}.1. 
By Lemma \ref{finfact}.3, the push-forward $\overset\gets g_\ast$ preserves 
constructibility. 

The 
 finiteness of tor-dimension follows from Proposition \ref{pr}.1 and Lemma \ref{finfact}.2. 

Let $\mathcal{G}:=\overset\gets g_\ast(R\Psi_f(K)
)$. 
By (\ref{stalkpush}) and the fact that $g$ is finite, for a point $(z'\leftarrow t)$ of $
Z'\overset\gets\times_SS$, we have a canonical isomorphism 
\begin{equation*}
\mathcal{G}_{(z'\gets t)}\cong 
\bigoplus_{x\in Z\times_{Z'}z'}R\Psi_f(K)_{(
x\leftarrow t)}. 
\end{equation*}
We show that the restriction $\mathcal{G}|_{Z'\overset\gets\times_S(S\setminus Z')}=
\overset\gets g_\ast(R\Psi_f(K)
|_{Z\overset\gets\times_S(S\setminus Z')})$ is locally constant. To prove this, it is 
 enough to show that, for points $(z'_1\gets t_1)$, $(z_2'\gets t_2)$ of $Z'
 \overset\gets\times_S(S\setminus Z')$ and a 
 specialization $(z'_1\gets t_1)\gets(z'_2\gets t_2)$, the morphism of stalks 
 $\mathcal{G}_{(z'_1\gets t_1)}\to \mathcal{G}_{(z'_2\gets t_2)}$ is an isomorphism by Lemma \ref{noe} below since 
the topos $Z'\overset\gets\times_S(S\setminus Z')$ is noetherian 
(\cite[Lemme 9.3.]{Org}). 
We prove this by showing that 
$\mathcal{G}_{(z'_1\gets t_1)}\to
\mathcal{G}_{(z'_1\gets t_2)}$ and $\mathcal{G}_{(z'
_1\gets t_2)}\to
\mathcal{G}_{(z'_2\gets t_2)}$ are both isomorphisms. 

By Proposition \ref{pr}.2, we have a distinguished triangle 
\begin{equation*}
\mathcal{G}_{(z'_1\gets t_1)}\to
\mathcal{G}_{(z'_1\gets t_2)}\to
\bigoplus_{x\in Z\times_{Z'}z'_1}
\bigoplus_{z\in Z_{(x)}\times_{S_{(z'_1)}}t_1}
R\Phi_f(K)_{(z\gets t_2)}
\to.
\end{equation*}
 Since $t_1$ is a geometric point of $S\setminus Z'$, $
Z_{(x)}\times_{S_{(z'_1)}}t_1$ is empty. 
This shows the former isomorphism 
 $\mathcal{G}_{(z'_1\gets t_1)}\overset\cong\to
\mathcal{G}_{(z'_1\gets t_2)}$. 

To prove the latter isomorphism, also by Proposition \ref{pr}.2, we have a 
distinguished triangle
\begin{equation}\label{d}
\mathcal{G}_{(z'_1\gets z'_2)}\to
\mathcal{G}_{(z'_1\gets t_2)}\to\bigoplus_
{x\in Z\times_{Z'}z'_1}\bigoplus_{z\in Z_{(x)}\times_{S_{(z'_1)}}z'_2}
R\Phi_f(K)
_{(z\gets t_2)}
\to.
\end{equation}
Since the restriction $K\lvert_{f^{-1}(Z')}$ is acyclic, the third term is naturally identified with 
\begin{equation*}
\bigoplus_{z\in Z\times_{Z'}z'_2}
R\Phi_f(K)
_{(z\gets t_2)}\cong
\bigoplus_{z\in Z\times_{Z'}z'_2}
R\Psi_f(K)
_{(z\gets t_2)}\cong
\mathcal{G}_{(z'_2\gets t_2)}. 
\end{equation*}
Note that the restriction of $\mathcal{G}$ to $Z'\overset\gets\times_{Z'}Z'\subset Z'\overset\gets\times_SS$ is acyclic 
since $K\lvert_{f^{-1}(Z')}$ is acyclic 
and the formations of $\overset\gets g_\ast$ and $R\Psi_f(K)$ commute with the base change $Z'\to S$. 
In particular, the first term of (\ref{d}) is acyclic. 
The assertion follows. 
\qed
}
\vspace{0.1in}

To give a proof of Lemma \ref{noe}, we need to recall some facts on 
noetherian topoi. Let $T$ be a noetherian topos. Recall that a point of 
$T$ is a morphism $t\colon ({\rm Set})\to T$ of topoi from the category of sets. 
For an object $X$ of $T$ and a point $t$ of $T$, denote $X_t:=t^\ast X$. 
For points $s$ and $t$, we call a natural transformation $s^\ast\to t^\ast$ a 
specialization $t\to s$. Let $|T|$ be the set of isomorphism classes of 
points of $T$. We denote the final object of $T$ by the same letter $T$. 
For a subobject $U$ of $T$, denote by $|U|$ the set of points $t\in|T|$ 
such that $U_t$ is not empty. There is a topology on $|T|$ whose open 
subsets are precisely $|U|$ for subobjects $U$ of $T$. This makes $|T|$ 
a noetherian topological space. The assignment $U\mapsto|U|$ gives a 
bijection between the set of subobjects of $T$ and  the set of open subsets of 
$|T|$, since $T$ has enough points (\cite[Proposition 9.0]{point}). 

Let $t$ be a point of $T$. The stalk $X_t$ of an object $X$ of $T$ at $t$ 
can be computed as follows (\cite[6.8]{SGA4}).  
Let ${\rm Nbd}(t)$ be the category whose objects are pairs $(U,x)$ where 
$U$ are quasi-compact objects of $T$ and $x$ are elements of $U_t$. The morphisms 
$(U_1,x_1)\to(U_2,x_2)$ are morphisms $U_1\to U_2$ which send 
$x_1$ to $x_2$. There is a functorial isomorphism 
\begin{equation*}
\varinjlim_{(U,x)\in {\rm Nbd}(t)}{\rm Hom}_T(U,X)\cong X_t. 
\end{equation*}
\begin{lm}\label{noe}
Let $T$ be a noetherian topos. Assume that, for a quasi-compact 
object $U$ of T and a point $t$ of $T$, the stalk $U_t$ is finite. 
\begin{enumerate}
\item For points $s$ and $t$ of $T$, the following are equivalent. 
\begin{enumerate}
\item There is a specialization $t\to s$. 
\item $s\in\overline{\{t\}}\subset |T|$. 
\end{enumerate}
\item Let $X$ be an object of 
$T$ whose stalks at all points are finite sets. The following 
are equivalent. 
\begin{enumerate}
\item $X$ is locally constant constructible. 
\item For all specializations $t\to s$ of points of $T$, 
the canonical maps $X_s\to X_t$ are bijective. 
\end{enumerate}
\end{enumerate}
\end{lm}
\proof{
Although these are well-known, we include proofs since we cannot find a 
reference. 

1. The implication $(a)\Rightarrow (b)$ is obvious from the 
definition of the topology of $|T|$. 

We prove $(b)\Rightarrow (a)$. For all objects $(U,x)$ of ${\rm Nbd}(s)$, 
$U_t$ are non-empty finite sets. Hence $\varprojlim_{(U,x)\in {\rm Nbd}(s)}
U_t$ is non-empty since ${\rm Nbd}(s)^{\rm op}$ is filtered. An element of 
$\varprojlim_{(U,x)\in {\rm Nbd}(s)}U_t$ gives a functor 
${\rm Nbd}(s)\to {\rm Nbd}(t)$, which defines a specialization $t\to s$. 

2. First we prove $(a)\Rightarrow(b)$. Take a covering $(T_i\to T)_i$ of 
$T$ so that $X\times_TT_i$ is a disjoint union of $T_i$. 
Let $t\to s$ be a specialization of points. Take $i$ so that $T_{i,s}$ is non-
empty. A lift $\bar{s}$ of $s$ to a point of $T_i$ gives a lift of $t$ to a point 
$\bar{t}$ of $T_i$ and a specialization $\bar{t}\to\bar{s}$. 
Since the stalk $(X\times_TT_i)_{\bar{s}}$ (resp. $(X\times_TT_i)_{\bar{t}}$) 
is canonically isomorphic to $X_s$ (resp. $X_t$), the assertion follows. 

Next we show $(b)\Rightarrow(a)$. Fix a point $t$ of $T$. 
Take a quasi-compact object $U$ of $T$ and $x\in U_t$ so that 
the natural map ${\rm Hom}_T(U,X)\to X_t$ is surjective. 
Replacing $(T,t)$ by $(U,x)$, we may assume that there is a morphism 
$\coprod T\to X$, from the disjoint union of finitely many copies of $T$, 
which gives a bijection $\coprod T_t\to X_t$. 
Further replacing $T$ by the connected component which contains $t$, 
we may assume that $|T|$ is connected. Then it follows that 
the morphism $\coprod T\to X$ is an isomorphism from 
1 and the fact that $T$ has enough points (\cite[Proposition 9.0]{point}). 
\qed
}

\subsection{Tame symbols for \'etale sheaves}

\begin{df}\label{tameobj1}
Let $S$ be a noetherian scheme. 
Let $C$ be a separated smooth $S$-curve. Let $Z\subset C$ be a 
closed subscheme finite \'etale over $S$. 
For a point $z$ of $Z$, denote by $\eta_z$ 
the generic point of the henselization of $C_{k(s),(z)}$ where $s\in S$ is 
the image of $z$. 
For a locally constant constructible object 
$X$ of $Z\overset \gets\times_{C}(C\setminus Z)$, 
we say that $X$ is a {\rm tame object} if, for every point $z\in Z$, 
the restriction of $X$ to $\eta_z\cong z\overset \gets\times_{C_{k(s)}}
(C_{k(s)}\setminus z)\subset Z\overset \gets\times_{C}(C\setminus Z)$ is tamely ramified. 
\end{df}
By \cite[Proposition 5.5]{Aby}, a locally constant constructible object 
$X$ is a tame object if and only if, for every $generic$ point $z\in Z$, 
the restriction to $\eta_z$ is tamely ramified. 

Let $f\in\Gamma(C,\mathcal{O}_C)$ be a global section which 
generates the ideal sheaf of $Z$. 
We construct, from a tame object of $Z\overset\gets\times_{C}
(C\setminus Z)$ and such a section $f$, a locally constant 
constructible object 
$\langle X,f\rangle$ of the \'etale topos of $Z$ when $S$ is normal. 
We start with more general setting \cite[1.7.8.]{Weil2}. 
Let $Y$ be a regular scheme (resp. a smooth scheme over a scheme $S$). Let 
$D\subset Y$ be a regular divisor (resp. a smooth divisor relative to $S$). 
Denote the complement by $U$. 
Denote by $i\colon D\to Y,j\colon U\to Y$ the immersions. 
Assume that the ideal sheaf of $D$ is globally generated by $z\in 
\Gamma(Y,\mathcal{O}_Y)$. Let $\mathcal{F}$ be a locally constant 
constructible sheaf of sets on $U$ which is tamely ramified 
along $D$. For an integer $n\geq1$ which is invertible in $Y$, 
let $Y_n:={\rm Spec}(\mathcal{O}_Y[t]/(t^n+z))$. This is a finite 
totally tamely ramified 
covering of $Y$ with a lifting $D\to Y_n$ of $D\to Y$. Denote by  $U_n:=U\times_YY_n$ the complement of $D$ in $Y_n$. 
Zariski locally on $Y$, we can find such  an $n$ that $\mathcal{F}|_{U_n}$ extends to a locally constant sheaf on $Y_n$, which we denote by $\mathcal{F}_n$. 
\begin{df-lm}\label{tameext}
Let the notation be as above. 
\begin{enumerate}
\item The restriction of $\mathcal{F}_n$ to $D\subset Y_n$ is 
independent of the choice of $n$. 
Hence the restrictions glue to a locally constant constructible sheaf 
on $D$, which we denote by $\langle\mathcal{F},z\rangle$. 
\item Suppose that $Y$ and $D$ are smooth over a scheme $S$. 
Then, the formation of $\langle\mathcal{F},z\rangle$ commutes 
with base change $S'\to S$. 
\end{enumerate}
\end{df-lm}
\proof{

1. Let $n,m\geq1$ be integers. The assertion follows from 
the fact that there is an $\mathcal{O}_Y$-morphism 
$\mathcal{O}_Y[t]/(t^n+z)\to\mathcal{O}_Y[u]/(u^{nm}+z)$ 
sending $t\mapsto u^m$. 

2. It follows since the fibered product $Y_n\times_SS'$ is canonically 
isomorphic to 
$(Y\times_SS')_n$. 
\qed
}

We go back to the situation of vanishing topoi. Let $S$ be a 
noetherian normal scheme. Let $f\in\Gamma(C,\mathcal{O}_C)$ be a global section which 
generates the ideal sheaf of $Z$ which is finite \'etale over $S$. 
Let $X$ be a tame object on 
$Z\overset\gets\times_{C}
(C\setminus Z)$. 
For a point $z\in Z$, the restriction of $X$ to 
$C_{(z)}\setminus Z_{(z)}\cong z\overset\gets\times
_{C}(C\setminus Z)\subset Z\overset\gets\times_{C}
(C\setminus Z)$ is a locally constant constructible 
sheaf tamely ramified along $Z_{(z)}$. 
\begin{df-lm}\label{tameobj}
Assume that $S$ is noetherian normal.
\begin{enumerate}
\item  
Locally constant constructible sheaves $\langle X|_{C_{(z)}
\setminus Z_{(z)}},
f\rangle$ on $Z_{(z)}$ glue to a locally constant constructible sheaf on 
$Z$. 
We denote this sheaf by $\langle X,f\rangle$. 
\item The formation of $\langle X,f\rangle$ commutes 
with arbitrary base change $S'\to S$. 
\end{enumerate}
\end{df-lm}
\proof{
1. 
Let $t\in Z$ be a generic point and $z\in Z$ be a point which is a specialization 
of $t$.  
The restriction of $\langle X|_{C_{(t)}\setminus t},f\rangle$ to $Z_{(z)}\times_Zt$ and that of $\langle X|_{C_{(z)}\setminus Z_{(z)}},
f\rangle$ to $Z_{(z)}\times_Zt$ coincide by Lemma \ref{tameext}.2. 
The assertion follows. 

2. It follows from Lemma \ref{tameext}.2. 
\qed
}

\subsection{Flat functions and trace maps}

We give a definition of flat functions. 
\begin{df}\label{flmoi}(\cite[Definition 2.1]{Sai17})
Let $h\colon Z\to S$ be a quasi-finite morphism of schemes. A function $\varphi\colon Z\to \mathbb{Z}$ is said to be 
{\rm flat} if, for every geometric point $(x\gets t)$ of $Z\overset\gets\times_SS$, 
we have 
\begin{equation*}
\varphi(x)=\sum_{z\in Z_{(x)}\times_{S_{(s)}}t}\varphi(z), 
\end{equation*}
where $s$ is the image $h(x)$. 
For a geometric point $z$ of $Z$, we write 
$\varphi(z)$ for the value of $\varphi$ at the image of $z\to Z$. 
\end{df}

\begin{ex}\label{flafcn}
Let $Z$ and $S$ be as in Definition \ref{flmoi}. 
\begin{enumerate}
\item Let $i\colon S'\to S$ be a morphism of schemes and $
i_Z\colon Z'=Z\times_SS'\to Z$ be the projection. 
If the function $\varphi\colon Z\to \mathbb{Z}$ is flat over $S$, the composition $\varphi\circ i_Z$ is flat over $S'$. 
\item Assume that $\mathcal{O}_Z$ is of finite tor-dimension as 
an $h^{-1}\mathcal{O}_S$-module. 
For a point $z\in Z$, write $\bar{z}$ for a geometric point above $z$ and 
$\bar{s}$ for the geometric point of $S$ defined by $\bar{z}\to Z\to S$.  
Then the function 
$Z\to\mathbb{Z}$ defined by $z\mapsto\sum_i(-1)^i{\rm dim}_{k(\bar{s})}
{\rm Tor}_i^{\mathcal{O}_{S,(\bar{s})}}(\mathcal{O}_{Z,(\bar{z})},k(\bar{s}))$ 
is flat over $S$. 
\end{enumerate}
\end{ex}

\begin{lm}\label{tracefl}
Let $h\colon Z\to S$ be a separated quasi-finite morphism of 
noetherian schemes. Let $\varphi\colon Z\to \mathbb{Z}$ be a flat function over $S$. 
\begin{enumerate}
\item There exists a unique morphism ${\rm Tr}_\varphi\colon h_!\mathbb{Z}\to
\mathbb{Z}$ of \'etale sheaves on $S$ such that, for every geometric point $s$ 
of $S$, the stalk ${\rm Tr}_{\varphi,s}\colon\bigoplus_{z\in Z_{s}}\mathbb{Z}\to \mathbb{Z}$ 
sends $(n_z)_z$ to $\sum_zn_z\varphi(z)$. 
\item Let $i\colon S'\to S$ be a morphism of noetherian schemes and let 
$Z\xleftarrow{i_Z}Z'=Z\times_SS'\xrightarrow{h'}S'$ be the projections. 
Then the diagram 
\begin{equation*}
\xymatrix{
i^\ast h_!\mathbb{Z}\ar[rr]^{i^\ast {\rm Tr}_\varphi}\ar[d]_{\cong}&&
i^\ast \mathbb{Z}\ar[d]_{\cong}\\h'_!i_Z^\ast\mathbb{Z}\ar[rr]^{{\rm Tr}_{\varphi\circ i_Z}}
&&\mathbb{Z}, 
}
\end{equation*}
where the vertical arrows are canonical isomorphisms, is commutative. 
\end{enumerate}
\end{lm}
\proof{
They are done in \cite[Proposition 6.2.5]{SGA4-3}. 
\qed
}
\vspace{0.1in}

Let $S$ be a noetherian $\F_p$-scheme. 
Let $\Lambda$ be a finite local ring in which $p$ is invertible. The Artin-Schreier covering 
$\mathbb{A}^1_S\to\mathbb{A}^1_S$ defined by $x\mapsto x^p-x$ is a Galois covering 
whose Galois group is canonically isomorphic to $\F_p$. For a non-trivial character 
$\psi\colon\F_p\to\Lambda^\times$, this covering defines a locally constant sheaf $\mathcal{L}_\psi(x)$ of 
invertible $\Lambda$-modules on $\mathbb{A}^1_S$, which is the Artin-Schreier sheaf. 
For a section $S\to\mathbb{A}^1_S$, we denote by $\mathcal{L}_\psi(f\cdot x)$ the pull-back of the Artin-Schreier 
sheaf by the  multiplication-by-$f$ map $\mathbb{A}^1_S\to\mathbb{A}^1_S$. 
\begin{lm}\label{AStr}
Let the notation be as above. Fix a non-trivial character $\psi\colon \F_p\to\Lambda^\times$. 
Let $Z$ be a finite $S$-scheme. 
For a $Z$-morphism 
$f\colon Z\to\mathbb{A}^1_Z$ and a flat function $\varphi\colon Z\to \mathbb{Z}$ over $S$, 
there exists an element $\mathcal{L}_\psi(\varphi\cdot f)\in{\rm H}^1(\mathbb{A}^1_S,\Lambda^\times)$  
with the following properties. 
\begin{enumerate}
\item For a morphism $i\colon S'\to S$ of noetherian schemes, $\mathcal{L}_\psi(\varphi\circ i_Z\cdot f')$ 
coincides with the image of 
$\mathcal{L}_\psi(\varphi\cdot f)$ by ${\rm H}^1(\mathbb{A}^1_S,\Lambda^\times)\to {\rm H}^1(\mathbb{A}^1_{S'},\Lambda^\times)$. 
Here $i_Z\colon Z':=Z\times_SS'\to Z$ is the projection and 
$f'\colon Z'\to\mathbb{A}^1_{Z'}$ is the base change of $f$. 
\item When $S={\rm Spec}(k)$ is the spectrum of a perfect field $k$, 
$\mathcal{L}_\psi(\varphi\cdot f)$ is equal to $\mathcal{L}_\psi((\sum_{z\in Z}\varphi(z){\rm Tr}_{k(z)/k}f(z))\cdot x)$, 
where $x$ is the standard coordinate of $\mathbb{A}^1_k$. 
\end{enumerate}
\end{lm}
\proof{
Let $\varphi\circ p_Z\colon \mathbb{A}^1_Z\to Z\to\mathbb{Z}$ be the composition of 
$\varphi$ and the projection  $p_Z\colon \mathbb{A}^1_Z\to Z$. This is a flat function over $\mathbb{A}^1_S$. 
Let ${\rm Tr}_{\varphi\circ p_Z}\colon (h\times{\rm id})_\ast\mathbb{Z}\to\mathbb{Z}$ be the trace map constructed in Lemma \ref{tracefl} 
from $\varphi\circ p_Z$ and the base change $h\times{\rm id}\colon\mathbb{A}^1_Z\to\mathbb{A}^1_S$ of $h$. 
Tensoring $\Lambda^\times$ and taking ${\rm H}^1$, we get a group homomorphism 
\begin{equation}\label{trmapl}
{\rm H}^1(\mathbb{A}^1_Z,\Lambda^\times)\to{\rm H}^1(\mathbb{A}^1_S,\Lambda^\times). 
\end{equation}
Let $\mathcal{L}_\psi(f\cdot x)\in{\rm H}^1(\mathbb{A}^1_Z,\Lambda^\times)$ 
be the pull-back of the rank $1$ locally constant $\Lambda$-sheaf $\mathcal{L}_\psi(x)$ on $\mathbb{A}^1_Z$ 
by the map $\mathbb{A}^1_Z\to\mathbb{A}^1_Z$ 
defined by $x\mapsto f\cdot x$. Define $\mathcal{L}_\psi(\varphi\cdot f)$ to be the image 
of $\mathcal{L}_\psi(f\cdot x)$ by (\ref{trmapl}). 
By Lemma \ref{tracefl}.2, this satisfies the properties. 
\qed
}

\section{Local Epsilon Factors (cf. \cite{Del}, \cite{Lau}, \cite{Y3})}\label{LEF}

In this preliminary section, we review theories of local epsilon factors for 
henselian traits of equal-characteristic. 
We fix prime numbers $p$ and $\ell$ so that $p\neq\ell$. 
Let $\Lambda$ be a finite local ring with residual characteristic $\ell$.  
We also  fix a non-trivial character $\psi\colon \mathbb{F}_p\to \Lambda^\times$.
\subsection{Generalities on local epsilon factors}
Let $k$ be a perfect field of characteristic $p$. 
Let $T$ be a henselian trait which is isomorphic to the henselization of $\mathbb{A}^1_{k}$ at a 
closed point. Let $s$ and $\eta$ be the closed point and the 
generic point of $T$ respectively. When $k$ is a finite field, in \cite{Del}, a constant $\varepsilon_{\psi}(T,\mathcal{F},\omega)\in E^\times$, 
or simply $\varepsilon(T,\mathcal{F},\omega)$, is defined  for a constructible complex $\mathcal{F}\in D^b_c(T,
E)$ and a non-zero rational $1$-form $\omega\in 
\Omega_{k(\eta)}^1$. Here $E$ is a finite extension of $\mathbb{Q}_\ell$ and $\psi\colon\F_p\to E^\times$ is a fixed non-trivial 
character.

\begin{thm}(\cite{Del}, \cite{Lau})\label{ep}
 Let $T$ be as above. Assume that $k$ is finite. 
Fix a non-trivial character $\psi\colon\mathbb{F}_p\to 
E^\times$. For each complex $\mathcal{F}\in D^b_c(T,E)$ and a non-zero 
rational $1$-form $\omega\in \Omega^1_{k(\eta)}$, we can attach, in a 
canonical way, an element 
\begin{equation*}
\varepsilon(T,\mathcal{F},\omega)\in E^\times
\end{equation*}
which satisfies the following properties: 
\begin{enumerate}
\item The element $\varepsilon(T,\mathcal{F},\omega)$ only depends on the 
isomorphism class of $(T,\mathcal{F},\omega)$. 
\item For an distinguished triangle
\begin{equation*}
\mathcal{F}_1\to \mathcal{F}_2\to \mathcal{F}_3\to
\end{equation*}
in $D^b_c(T,E)$, we have $\varepsilon(T,\mathcal{F}_2,\omega)=
\varepsilon(T,\mathcal{F}_1,\omega)\cdot\varepsilon(T,\mathcal{F}_3,\omega)$.
\item If $\mathcal{F}$ is supported on the closed point $s$, we have 
\begin{equation*}
\varepsilon(T,\mathcal{F},\omega)={\rm det}(-{\rm Frob}_s,\mathcal{F}_{\bar{s}})^{-1}.
\end{equation*}
\item Let $\eta_1/\eta$ be a finite separable extension and $f\colon 
T_1\to T$ be the 
normalization in $\eta_1$. For a constructible complex $\mathcal{F}_1\in D^b_c(T_1,E)$ with 
generic rank $0$, i.e. $\rk \mathcal{F}_{1,\bar{\eta_1}}=0$, we have 
\begin{equation*}
\varepsilon(T_1,\mathcal{F}_1,f^\ast\omega)=\varepsilon(T,f_\ast \mathcal{F}_1,\omega).
\end{equation*}
\item Let $\mathcal{G}$ be a smooth $E$-sheaf on $\eta$ of rank $1$, which  induces a character $\chi\colon k(\eta)^\times\to E^\times$ via 
the Artin map $k(\eta)^\times\to G_\eta^{ab}$. Let $j\colon\eta\to T$ 
be the immersion. Then, we have 
\begin{equation*}
\varepsilon(T,j_\ast\mathcal{G},\omega)=\varepsilon(\chi,\Psi_\omega). 
\end{equation*}
Here $\Psi_\omega(a)=\psi\circ{\rm {\rm Tr}}_{k(s)/\mathbb{F}_p}({\rm Res}(a\cdot
\omega))$ and $\varepsilon(\chi,\Psi_\omega)$ is the Tate constant 
(\cite[3.1.3.2]{Lau}). 
\end{enumerate}
\end{thm}
Here we fix a normalization 
  of the local class field theory so that the Artin map sends 
  a uniformizer to a geometric Frobenius.
We identify a smooth $E$-sheaf on $\eta$ and a finite dimensional 
$E$-representation of the absolute Galois group $G_\eta$. 
For a smooth $E$-sheaf $V$ on $\eta$, we also denote by $\varepsilon_0(T,V,\omega)$ the constant  
$\varepsilon(T,V_!,\omega)$, where $V_!$ is the $0$-extension of $V$ to $T$. 

The local epsilon factors admit the following properties. 
\begin{pr}\label{elemep}
\begin{enumerate}
\item (\cite[3.1.5.5]{Lau}) For a non-zero element $a\in\Gamma(\eta,\mathcal{O}_\eta)$, we have 
\begin{equation*}
\varepsilon(T,\mathcal{F},a\cdot\omega)=
{\rm det}(\mathcal{F}_\eta)(a)\| a\|^{-{\rm rk}\mathcal{F}_\eta}
\varepsilon(T,\mathcal{F},\omega).
\end{equation*}
\item (\cite[3.1.5.6]{Lau}) Let $\mathcal{G}$ be a smooth $E$-sheaf on $T$. Then, we have 
\begin{equation*}
\varepsilon(T,\mathcal{F}\otimes\mathcal{G},\omega)=
{\rm det}(\mathcal{G})({\rm Frob}_s)^{a(T,\mathcal{F},\omega)}\cdot
\varepsilon(T,\mathcal{F},\omega)^{{\rm rk}\mathcal{G}}. 
\end{equation*}
Here $a(T,\mathcal{F},\omega)$ is defined to be ${\rm rk}(\mathcal{F}_\eta)+
{\rm Sw}(\mathcal{F}_\eta)-{\rm rk}(\mathcal{F}_s)+
{\rm rk}(\mathcal{F}_\eta)\cdot{\rm ord}(\omega)$. \label{unrtw}
\end{enumerate}
\end{pr}
Note that, with a constructible complex $\mathcal{F}\in D^b_c(T,\mathcal{O}_E)$, we can define  
a constant $\varepsilon(T,\mathcal{F},\omega)$ to be $\varepsilon(T,\mathcal{F}\otimes_{\mathcal{O}_E}E,\omega)$. 

Yasuda gives a generalization of this theory \cite{Y1}, \cite{Y3}. 
In his setting, we can take a finite local ring as the coefficient ring of sheaves and can 
treat a  henselian trait with a general perfect residue field (which further can be of mixed characteristic). 

We adjust his result to our setting. 
Let $\Lambda$ be a finite local ring of residual characteristic $\ell$. 
Take and fix a non-trivial additive character $\psi\colon\mathbb{F}_p\to\Lambda^\times$. 
Let $k$ be a perfect 
field of characteristic $p$ and $T$ be a henselian trait isomorphic to the henselization of $\mathbb{A}^1_k$ at a closed point. We denote its generic point by $\eta$. We denote the absolute Galois group of $k$ (resp. $\eta$) 
by $G_k$ (resp. $G_\eta$). 
\begin{thm}(\cite[4.12]{Y3})\label{torloc}
Let the notation be as above. For a triple $(T,(\rho,V),\omega)$ where $V$is a finite free $\Lambda$-module with a continuous 
group homomorphism $\rho\colon G_\eta\to{\rm GL}(V)$ and $\omega\in\Omega^1_\eta$ is a non-zero rational section, there is a canonical way to attach a 
continuous character $\varepsilon_{0,\Lambda}(T,V,\omega)\colon G_k^{ab}\to\Lambda^\times$ with the following properties. 
\begin{enumerate}
\item The character only depends on the isomorphism class of $(T,(\rho,V),\omega)$. 
\item For a short exact sequence $0\to V'\to V\to V''\to0$ of representations of $G_\eta$, we have 
\begin{equation*}
\varepsilon_{0,\Lambda}(T,V,\omega)=\varepsilon_{0,\Lambda}(T,V',\omega)\cdot\varepsilon_{0,\Lambda}(T,V'',\omega). 
\end{equation*}
\item For a local ring homomorphism $f\colon \Lambda\to\Lambda'$, we have 
\begin{equation*}
f\circ\varepsilon_{0,\Lambda}(T,V,\omega)=\varepsilon_{0,\Lambda'}(T,V\otimes_\Lambda\Lambda',\omega)
\end{equation*}
as characters $G_k\to\Lambda'^\times$. 
\item Assume that the residue field $k$ of $T$ is finite and that there exists a local ring morphism 
$f\colon \mathcal{O}_E\to\Lambda$ from the ring of integers of a finite extension $E/\mathbb{Q}_\ell$ 
such that $V$ comes from a representation on $\mathcal{O}_E$, i.e. there is a representation $V'$ of $G_\eta$ 
on a finite free $\mathcal{O}_E$-module such that $V'\otimes_{\mathcal{O}_E}\Lambda\cong V$. 
Then we have 
\begin{equation*}
\varepsilon_{0,\Lambda}(T,V,\omega)({\rm Frob_k})=(-1)^{\rk V+{\rm Sw}V}f(\varepsilon_{0}(T,V'\otimes_{\mathcal{O}_E}E,\omega)). 
\end{equation*}
Here the local epsilon factor in the right hand side is the one in Theorem \ref{ep}. 
\end{enumerate}
\end{thm}

\begin{rmk}
In the case that the residue field $k$ is finite, before the work of \cite{Y3}, 
Yasuda defines a local epsilon factor $\varepsilon_{0,\Lambda}(T,V,\omega)$ as an element of $\Lambda^\times$ in \cite{Y1}. 
The relation between them is as follows: 
\begin{equation*}
\varepsilon_{0,\Lambda}(T,V,\omega)({\rm Frob_k})=(-1)^{\rk V+{\rm Sw}V}\varepsilon_{0,\Lambda}(T,V,\omega). 
\end{equation*}
Here the local epsilon factor in the left hand side is given in Theorem \ref{torloc}. 
\end{rmk}

For more properties, see \cite{Y3}. By the multiplicativity, we also define a local epsilon factor 
$\varepsilon_{0,\Lambda}(T,K,\omega)$ for 
a complex $K\in D_{\rm ctf}(\eta,\Lambda)$. 

Let us explain precisely the relation between the local epsilon factors defined in \cite{Y3} and the one 
in Theorem \ref{torloc}. 

First let us explain a counterpart of an additive character in the classical setting. 
Let $F$ be the completion of the function field of $T$. Let $\mathfrak{m}_F$ be the maximal ideal 
of the ring of integers of $F$. 
For integers $m\leq n$, define a contravariant functor $F^{[m,n]}$ from the category of affine $k$-schemes to the 
category of sets to be ${\rm Spec}(A)\mapsto A\otimes_k\mathfrak{m}_F^m/\mathfrak{m}_F^{n+1}$. This functor is represented by an 
affine smooth group $k$-scheme. Indeed, fixing a uniformizer $\pi\in F$, it is isomorphic to the functor 
${\rm Spec}(A)\mapsto A^{\oplus m-n}$, which is represented by $\mathbb{A}^{m-n}_k$. 
We have a natural inclusion $F^{[m,n]}\to F^{[m-1,n]}$ defined by $\mathfrak{m}_F^m/\mathfrak{m}_F^{n+1}\to\mathfrak{m-1}_F^m/\mathfrak{m}_F^{n+1}$. {\it A non-trivial invertible character sheaf} 
$\tilde{\psi}^{[m,n]}$ on $F^{[m,n]}$ is a non-trivial invertible sheaf of $\Lambda$-modules on $F^{[m,n]}$ such that 
$m^\ast\tilde{\psi}^{[m,n]}$ is isomorphic to the external product $\tilde{\psi}^{[m,n]}\boxtimes\tilde{\psi}^{[m,n]}$, where 
$m\colon  F^{[m,n]}\times F^{[m,n]}\to F^{[m,n]}$ is the addition. Denote the set of isomorphism classes of non-trivial 
invertible character sheaves on $F^{[m,n]}$ by ${\rm ACh}^0(F^{[m,n]},\Lambda)$ \cite[4.1]{Y3}. 
The natural inclusion $F^{[m,n]}\to F^{[m-1,n]}$ induces a map ${\rm ACh}^0(F^{[m-1,n]},\Lambda)\to {\rm ACh}^0(F^{[m,n]},\Lambda)$ 
by the pull-back. {\it A non-trivial additive character sheaf} $\tilde{\psi}$ on $F$ is an element of 
\begin{equation*}
\coprod_{n\in\mathbb{Z}}\varprojlim_{m\leq -n-1}{\rm ACh}^0(F^{[m,-n-1]},\Lambda). 
\end{equation*}
When $\tilde{\psi}\in \varprojlim_{m\leq -n-1}{\rm ACh}^0(F^{[m,-n-1]},\Lambda)$, the integer $n$ is called the conductor of 
$\tilde{\psi}$. 

In \cite{Y3}, Yasuda defines a character $\tilde{\varepsilon}_{0,\Lambda}(V,\tilde{\psi})$ from a representation $(\rho,V)$ as 
above and a non-trivial additive character sheaf $\tilde{\psi}$ on $F$. For a non-zero rational $1$-form $\omega\in\Omega^1_\eta$, 
we denote the corresponding additive character sheaf by $\tilde{\psi}_\omega$. Then we define 
$\varepsilon_{0,\Lambda}(T,V,\omega)=\tilde{\varepsilon}_{0,\Lambda}(V,\tilde{\psi}_\omega)$. 

The construction of $\tilde{\psi}_\omega$ goes as follows. 
When the residue field $k$ is finite, $\omega$ defines an additive character $\psi_\omega\colon F\to\Lambda^\times$ by 
$a\mapsto\psi({\rm Tr}_{k/\mathbb{F}_p}\circ {\rm Res}(a\omega))$. Then $\tilde{\psi}_\omega$ is the 
additive character sheaf corresponding to $\psi_\omega$ \cite[Corollary 4.3]{Y3}. In general, take a uniformizer 
$\pi\in F$, which defines an inclusion $\mathbb{F}_p((\pi))\to F$. When  $\omega=d\pi$, define 
$\tilde{\psi}_\omega$ to be the pull-back of the additive character sheaf $\tilde{\psi}_{d\pi}$ on 
$\mathbb{F}_p((\pi))$ \cite[Corollary 4.4.2]{Y3}. When $\omega=ad\pi$ for a non-zero element $a\in F$, 
define $\tilde{\psi}_\omega$ to be the pull-back of $\tilde{\psi}_{d\pi}$ by the multiplication-by-$a$ map 
$F\to F$ \cite[4.1]{Y3}. 

The next aim of this section is to explain a formula of local epsilon factors analogous to Laumon's formula \cite[3.5.1.1]{Lau}. 

We slightly generalize Laumon's local Fourier transform to general base schemes. 
Let $S$ be a noetherian scheme over $\F_p$. Fix a finite local ring $\Lambda$ of residue characteristic $\ell$. For a non-trivial character $\psi\colon\mathbb{F}_p\to\Lambda^\times$, 
let $\mathcal{L}_\psi(x)$ be the Artin-Schreier sheaf on $\mathbb{A}^1_S$. Let 
$\mathcal{L}_\psi(x_1\cdot x_2)$ be the pull-back of $\mathcal{L}_\psi(x)$ by 
the multiplication $\mathbb{A}^1_{\F_p}\times_{\F_p}\mathbb{A}^1_{\F_p}\to \mathbb{A}^1_{\F_p}$, 
where $x_1$ and $x_2$ are the standard coordinates of the first and second 
factors respectively.
We define $\bar{\mathcal{L}}_{\psi}$ to be the $0$-extension 
of $\mathcal{L}_\psi(x_1\cdot x_2)$ to $\mathbb{A}^1_{\F_p}\times_{\F_p}\mathbb{P}^1_{\F_p}$. We also use the same notation $\bar{\mathcal{L}}_{\psi}$ for 
the pull-back by $\mathbb{A}^1_S\times_S\mathbb{P}^1_S\to\mathbb{A}^1_{\F_p}\times_{\F_p}\mathbb{P}^1_{\F_p}$. 
Denote by $0_S$ (resp. $\infty_S$) the section $S\to\mathbb{A}^1_S$ (resp. 
$S\to\Proj^1_S$) at the origin (resp. at the infinity). Let 
\begin{equation*}
0_S\overset\gets\times_{\mathbb{A}^1_S}\mathbb{A}^1_S
\xleftarrow{\overset\gets p_1}(0_S\times_S\infty_S)\overset\gets\times_{\mathbb{A}^1_S\times_S\Proj^1_S}(\mathbb{A}^1_S\times_S\Proj^1_S)\xrightarrow{\overset\gets p_2}
\infty_S\overset\gets\times_{\Proj^1_S}\Proj^1_S
\end{equation*}
be the morphisms of topoi induced from the projections 
$\mathbb{A}^1_S\xleftarrow{p_1}\mathbb{A}^1_S\times_S\Proj^1_S\xrightarrow{p_2}\Proj^1_S$. 
Let $q\colon (0_S\times_S\infty_S)\overset\gets\times_{\mathbb{A}^1_S\times_S\Proj^1_S}(\mathbb{A}^1_S\times_S\Proj^1_S)\to\mathbb{A}^1_S\times_S\Proj^1_S$ be the second projection. 
\begin{df}\label{rellocfou}
Let $K\in D_{{\rm ctf}}(0_S\overset\gets\times_{\mathbb{A}^1_S}\mathbb{G}_{{\rm m},S},\Lambda)$ be a complex 
of finite tor-dimension whose cohomology sheaves 
are locally constant. 
Define the local Fourier transform $F^{(0,\infty)}(K)$ by 
\begin{equation*}
F^{(0,\infty)}(K):=R\overset\gets p_{2\ast}(\overset\gets p_1^\ast K_!\otimes^L_\Lambda q^\ast 
\bar{\mathcal{L}}_\psi)[1]|_{\infty_S\overset\gets\times_{\Proj^1_S}\mathbb{A}^1_S}\in D^b(\infty_S\overset\gets\times_{\Proj^1_S}\mathbb{A}^1_S,\Lambda). 
\end{equation*}
Here $K_!$ is the 
$0$-extension of $K$ by $0_S\overset\gets\times_{\mathbb{A}^1_S}\mathbb{G}_{{\rm m},S}\hookrightarrow0_S\overset\gets\times_{\mathbb{A}^1_S}\mathbb{A}^1_S$. 
\end{df}

\begin{pr}(cf. \cite[Proposition (2.4.2.2).]{Lau})\label{F}
Let the notation be as above. Let $\Lambda_0$ be the residue field of 
$\Lambda$. 
\begin{enumerate}
\item The formation of $F^{(0,\infty)}(K)$ commutes with arbitrary 
base change $S'\to S$. 
\item When $K$ has tor-amplitude in $[0,0]$, so is $F^{(0,\infty)}(K)$. 
\item\label{loccf} The local Fourier transform $F^{(0,\infty)}(K)$ is constructible 
and of finite tor-dimension. 
It is locally constant if and only if, 
for each $i\in\mathbb{Z}$, 
the function on $S$ defined by $s\mapsto{\rm dimtot}\mathcal{H}^i(K\otimes^L_{\Lambda}\Lambda_0)|_{\eta_{\bar{s}}}$, 
where $\bar{s}$ is the spectrum of an algebraic closure of $k(s)$ and 
$\eta_{\bar{s}}\cong0\overset\gets\times_{\mathbb{A}^1_{k(\bar{s})}}
\mathbb{G}_{{\rm m},k(\bar{s})}$ is the generic point of the 
henselization $\mathbb{A}^1_{k(\bar{s}),(0)}$ at the origin, 
is locally constant. 
\end{enumerate}
\end{pr}
\proof{
1. Let $s$ be a geometric point of $S$. 
Since the projection $\mathbb{A}^1_S\times_S\Proj^1_S\to\Proj^1_S$ 
is universally locally acyclic relatively to $\bar{\mathcal{L}}_\psi$ (\cite[Th\'eor\`eme (1.3.1.2)]{Lau}), so is 
the morphism 
$(\mathbb{A}^1_S\times_S\Proj^1_S)_{(0_s,\infty_s)}\to\Proj^1_{S,(\infty_s)}$ 
relatively to the restriction of $\bar{\mathcal{L}}_\psi$. 
Since $K$ is locally constant, the pair $(\mathbb{A}^1_S\times_S\Proj^1_S\to\Proj^1_S,0_S\times_S\infty_S
\xrightarrow{\cong}\infty_S)$ is 
universally locally acyclic relatively to $K_!\otimes^Lq^\ast\bar{\mathcal{L}}_\psi$ 
outside $0_S\times_S\Proj^1_S\subset\mathbb{A}^1_S\times_S\Proj^1_S$, in the sense of Definition \ref{relbcdef}.1. 
The assertion follows from Corollary \ref{relbc}, taking $X,S,Y,T,$ and $Z$ in the corollary as 
$\mathbb{A}^1\times_S\Proj^1_S,\Proj^1_S,0_S\times_S\infty_S,\infty_S,$ and $0_S\times_S\Proj^1_S$ respectively. 

2. Replacing $K$ by $K\otimes^L_{\Lambda}\Lambda_0$, 
we may assume that $\Lambda$ is a field. By 1, we may also assume that 
$S$ is the spectrum of a field. In this case, the assertion follows from the fact 
that 
the vanishing cycles functor preserves perversity up to shift \cite[Corollaire 4.6.]{Autor}. 

3. The complex $F^{(0,\infty)}(K)$ is of finite tor-dimension since the cohomological 
dimension of $\overset\gets p_{2\ast}$ is finite \cite[Proposition 3.1.]{Org}. 
For the constructibility and the local constancy, replacing $K$ by $K\otimes^L_{\Lambda}
\Lambda_0$, we may assume that $\Lambda$ is a field, and that $K$ is a locally constant sheaf. 
By 1 and the constructibility of $K$, we may assume that $S$ is of finite type 
over $\F_p$. 
Let $S'\to S$ be a proper surjective morphism of $\F_p$-schemes. 
Note that $F^{(0,\infty)}(K)$ is (resp. locally constant) constructible 
if and only if so is $F^{(0,\infty)}(K|_{S'})$. Indeed, the equivalence of 
the constructibility follows from 
\cite[Lemme 10.5.]{Org} and 1. By Lemma \ref{noe}, $F^{(0,\infty)}(K)$ 
is locally constant if and only if, for all specializations of points of 
$\infty_S\overset\gets\times_{\Proj^1_S}\mathbb{A}^1_S$ of types 
$(x_1\gets t)\gets(x_2\gets t)$ or $(x\gets t_1)\gets(x\gets t_2)$, 
the canonical maps of the stalks of $F^{(0,\infty)}(K)$ are isomorphisms. 
Hence the assertion follows from 1. 
Thus we may assume that $S$ is regular by de Jong's alteration theorem. 

Let $\eta\in S$ be a generic point. By the Gabber-Katz extension \cite{GK}, 
there exists a locally constant constructible sheaf $K'$ on $\mathbb{G}_{{\rm m},k(\eta)}$ 
which restricts to $K$ on $0\overset\gets\times_{\mathbb{A}^1_{k(\eta)}}
\mathbb{G}_{{\rm m},k(\eta)}$. By Zariski-Nagata's purity theorem, $K'$ extends to a locally constant constructible 
sheaf on 
$\mathbb{G}_{{\rm m},S}$, which is denoted also by $K'$. 
The projection $p_2\colon\mathbb{A}^1_S\times_S\Proj^1_S\to\Proj_S^1$ is universally locally acyclic relatively to $p_1^\ast K'_!\otimes^L\bar{\mathcal{L}}_\psi$ 
outside $0_S\times_S\Proj^1_S$ since so is $\bar{\mathcal{L}}_\psi$ 
(\cite[Th\'eor\`eme (1.3.1.2)]{Lau}) and 
$K'$ is locally constant. 
If the function $s\mapsto{\rm dimtot}K|_{\eta_{\bar{s}}}$ is locally constant, 
$\mathbb{A}^1_S\to S$ is 
universally locally acyclic relatively to $K'_!$ 
by \cite[Th\'eor\`eme 2.1.1]{semi-c}. Hence $\mathbb{A}^1_S\times_S\Proj^1_S\to \Proj^1_S$ is 
universally locally acyclic relatively to $p_1^\ast K'_!\otimes^L\bar{\mathcal{L}}_\psi$ outside $0_S\times_S\infty_S$. 
The assertion except the ``only if'' part follows from Propositions \ref{pr}.1, \ref{ula}. 

The ``only if'' part can be verified as follows. Assume that $F^{(0,\infty)}(K)$ is locally constant. 
Replacing $K$ by $K\otimes^L_\Lambda\Lambda_0$, we may assume that $\Lambda$ is a field. 
By 2, 
the function $\dt\mathcal{H}^i(K)|_{\eta_{\bar{s}}}$ coincides with 
${\rm rk}\mathcal{H}^iF^{(0,\infty)}(K)|_{\eta_{\bar{s}}}$ (Here we use Theorem \ref{lf}.2, which 
is proved independently of this proposition), hence the assertion. 
\qed
}

\vspace{0.1in}

Let $k$ be a perfect field of characteristic $p>0$, and 
 $T$ and $T^{\prime}$ be henselian traits of equal-characteristic 
with residue field $k$. Fix uniformizers $\pi$ and $\pi^{\prime}$ of $T$ and $T^{\prime}$ respectively.
Assume that the $k$-morphism $T\to\mathbb{A}^1_k$ 
(resp. $T'\to\Proj^1_k$) defined by $x\mapsto\pi$ (resp. $
x'\mapsto1/\pi'$), where $x$ (resp. $x'$) is the standard coordinate 
of $\mathbb{A}^1_k$ (resp. $\mathbb{A}^1_k\subset\Proj^1_k$), 
induces an isomorphism between $T$ (resp. $T'$) and the henselization 
of $\mathbb{A}^1_k$ (resp. $\Proj^1_k$) at $0$ (resp. at $\infty$). 
We identify $T$ and $T'$ with the henselizations. 
Let $s,\eta$ and $s',\eta^{\prime}$ be the closed points and 
the generic points of $T$ and $T^{\prime}$ respectively. 
We also denote the pull-back of $\bar{\mathcal{L}}_{\psi}$ by the map 
$T\times_kT'\to\mathbb{A}^1_k\times_k\mathbb{P}^1_k$ 
by the same letter $\bar{\mathcal{L}}_{\psi}$. 

Let $G_\eta$ and $G_{\eta'}$ be the absolute Galois groups of 
$\eta$ and $\eta'$. We identify the category $D_{\rm ctf}(\eta,\Lambda)$ (resp. $D_{\rm ctf}(\eta',\Lambda)$) with the derived category of 
complexes of continuous $G_\eta$-representations 
(resp. continuous $G_{\eta'}$-representations) on finite $\Lambda$-modules 
whose tor-dimensions are finite. 

There is a canonical upper numbering filtration $(G_{\eta'}^v)_v$ of 
$G_{\eta'}$ indexed by $v\in \mathbb{Q}_{\geq0}$. 
$G_{\eta'}^0$ is the inertia subgroup and $G_{\eta'}^{0+}:=\overline{
\cup_{v>0}G_{\eta'}^v}$ is the wild inertia subgroup. 
Let $D_{\rm ctf}(\eta',\Lambda)_{[0,1[}$ be the full subcategory 
of $D_{\rm ctf}(\eta',\Lambda)$ consisting of complexes on which 
$G_{\eta'}^1$ acts trivially. 

\begin{thm}(cf. \cite[Th\'eor\`eme (2.4.3)]{Lau})\label{lf} 
Let the notation be as above. 
\begin{enumerate}
\item $F^{(0,\infty)}$ induces an equivalence of categories 
\begin{equation*}
F^{(0,\infty)}\colon D_{\rm ctf}(\eta,\Lambda)\xrightarrow{\cong}
D_{\rm ctf}(\eta',\Lambda)_{[0,1[}.
\end{equation*}
\item For $K\in D_{\rm ctf}(\eta,\Lambda)$, we have 
\begin{itemize}
\item $\rk F^{(0,\infty)}(K)=\rk K+\sw K.$
\item $\sw F^{(0,\infty)}(K)=\sw K$. 
\end{itemize}
\end{enumerate}
\end{thm}
\proof{
Originally similar results for $\Ql$-sheaves are proven in \cite[Th\'eor\`eme (2.4.3)]{Lau}. Since the arguments in {\it loc. cit.} work, 
we only sketch the proof. 

1. Similarly as $F^{(0,\infty)}$, we define a functor $F^{(\infty,0)}\colon 
D_{\rm ctf}(\eta',\Lambda)\to D_{\rm ctf}(\eta,\Lambda)$ by 
$F^{(\infty,0)}(K):=R\overset\gets p_{1\ast}(\overset\gets p_2^\ast K_!\otimes^L_\Lambda
\bar{\mathcal{L}}_\psi)[1]|_\eta$, where we denote the projections by 
\begin{equation*}
s\overset\gets\times_TT\xleftarrow{\overset\gets p_1}
(s\times_ks')\overset\gets\times_{(T\times_kT')}(T\times_kT')
\xrightarrow{\overset\gets p_2}s'\overset\gets\times_{T'}T'. 
\end{equation*}
Similarly as \cite[Th\'eor\`eme (2.4.3)]{Lau}, we have natural 
isomorphisms $F^{(\infty,0)}\circ F^{(0,\infty)}(K)
\cong a^\ast K(-1)$ 
for $K\in D_{\rm ctf}(\eta,\Lambda)$, and 
$F^{(0,\infty)}\circ F^{(\infty,0)}(K)\cong a^\ast K(-1)$ for $K\in D_{\rm ctf}(\eta',\Lambda)_{[0,1[}$, where $a$ denotes the $k$-isomorphisms 
$T\to T$ and $T'\to T'$ which send $\pi$ and $\pi'$ to $-\pi$ and $-\pi'$ 
respectively. To show the assertion, then, it suffices to prove that the image of 
$F^{(0,\infty)}$ is contained in $D_{\rm ctf}(\eta',\Lambda)_{[0,1[}$. 
To this end, we may assume that $\Lambda$ is a field and that 
$K\in D_{\rm ctf}(\eta,\Lambda)$ is a sheaf. 
The local Fourier transform $V:=F^{(0,\infty)}(K)$ is a sheaf by Proposition \ref{F}.2. Similarly as \cite[Th\'eor\`eme (2.4.3)]{Lau}, we have an isomorphism $F^{(\infty,0)}(V^{G^1_{\eta'}})
\cong F^{(\infty,0)}(V)$. Hence the assertion follows from the isomorphisms 
\begin{equation*}
a^\ast V^{G^1_{\eta'}}(-1)\cong 
F^{(0,\infty)}\circ F^{(\infty,0)}(V)=F^{(0,\infty)}\circ F^{(\infty,0)}\circ 
F^{(0,\infty)}(K)\cong F^{(0,\infty)}(a^\ast K(-1))\cong a^\ast V(-1). 
\end{equation*}

2. The equalities follow from the existence of the Gabber-Katz extension and the 
Grothendieck-Ogg-Shafarevich formula. 
\qed
}

\begin{cor}\label{tame}
For a constructible complex $K\in D_{\rm ctf}(\eta,\Lambda)$, the determinant  
$\det(F^{(0,\infty)}(K))$ is tamely ramified. 
\end{cor}
\proof{
Since $\det(F^{(0,\infty)}(K))$ is of rank $1$ and 
$G^1_{\eta'}$ acts trivially on it, the assertion follows from Hasse-Arf Theorem. 
\qed
}

   \begin{thm}(cf. \cite[8.3]{Y2}, \cite[Th\'eor\`eme (3.5.1.1)]{Lau})\label{fou} 
   Let the notation be as in Theorem \ref{lf}. For a constructible complex $K\in D_{\rm ctf}(\eta,\Lambda)$, we have 
   \begin{equation}\label{locfouep}
   \varepsilon_{0,\Lambda}(T, K, d\pi ) =\langle \det(F^{(0, \infty)}(K)),\pi^{\prime}\rangle. 
   \end{equation}
\end{thm}
\proof{
We may assume that $K$ is represented by a locally constant sheaf $V$ of finite free $\Lambda$-modules. 
When the residue field $k$ is finite, it is proved in \cite[Proposition 8.3]{Y2}, and the argument given there 
is also valid in the general setting. 
\qed}

\begin{rmk}
In the sequel, it is harmless to take the formula (\ref{locfouep}) as the definition of 
local epsilon factors. 
\end{rmk}

To explain more properties on local epsilon factors by Yasuda, we need to recall the construction in \cite[4.2]{Y3}. 

Let $\chi\colon G_\eta\to M$ be a continuous group homomorphism into a finite abelian group  $M$. 
Fix a presentation 
\begin{equation}\label{repM}
0\to\mathbb{Z}^{\oplus n}\xrightarrow{\alpha}\mathbb{Z}^{\oplus n}\to M\to 0
\end{equation}
as an abelian group. This defines an injection ${\rm H}^1(\eta,M)\to{\rm H}^2(\eta,\mathbb{Z}^{\oplus n})$ of cohomology groups 
of $G_\eta$-modules with trivial action. 
Let $k(\eta)^{\rm sep}$ be a separable closure of $k(\eta)$ and $k(\eta)^{\rm ur}$ be the maximal 
unramified extension in it. Then the inflation morphism 
${\rm H}^2(k(\eta)^{\rm ur}/k(\eta),k(\eta)^{{\rm  ur}\times})\to {\rm H}^2(k(\eta)^{\rm sep}/k(\eta),k(\eta)^{{\rm  sep}\times})$ 
is an isomorphism. Thus we have a sequence of maps 
\begin{align}\label{paring}
k(\eta)^\times\times{\rm H}^1(\eta,M)\to k(\eta)^\times\times{\rm H}^2(\eta,\mathbb{Z}^{\oplus n})\to{\rm H}&^2(k(\eta)^{{\rm sep}}/
k(\eta),k(\eta)^{{\rm sep}\times})^{\oplus n}\\\notag
&\cong{\rm H}^2(k(\eta)^{\rm ur}/k(\eta),k(\eta)^{{\rm  ur}\times})^{\oplus n}\to{\rm H}^2(k,\mathbb{Z}^{\oplus n}). 
\end{align}
Here the second one is the cup product and the last one is induced from the valuation. 
Since ${\rm H}^1(\eta,M)$ (resp. ${\rm H}^1(k,M)$) is the kernel of $\alpha\colon {\rm H}^2(\eta,\mathbb{Z}^{\oplus n})\to
{\rm H}^2(\eta,\mathbb{Z}^{\oplus n})$ (resp. ${\rm H}^2(k,\mathbb{Z}^{\oplus n})\to{\rm H}^2(k,\mathbb{Z}^{\oplus n})$), 
The image of (\ref{paring}) lands into ${\rm H}^1(k,M)$. In this way, we have a pairing 
$k(\eta)^\times\times{\rm H}^1(\eta,M)\to{\rm H}^1(k,M)$. This is independent of the choice of (\ref{repM}) and 
natural in $M$. 
\begin{df}
For a pair $(a,\chi)\in k(\eta)^\times\times{\rm H}^1(\eta,M)$, we denote by 
$\chi_{[a]}\in{\rm H}^1(k,M)$ the image by this pairing. It is a character $G_k\to M$. 
\end{df}
\begin{lm}\label{compar}
Let $T$ be a henselian trait of equal-characteristic $p>0$ with 
perfect residue field $k$. 
Let $\chi\colon G_\eta\to M$ be a continuous group homomorphism into a finite abelian group $M$ and 
denote by $\mathcal{F}$ the corresponding locally constant sheaf on $\eta$. 
Assume that $\chi$ (hence also $\mathcal{F}$) is tamely ramified. 
Then, for a uniformizer $z\in k(\eta)^\times$, the character $\chi_{[z]}$ corresponds to 
$\langle\mathcal{F},z\rangle$ in Definition \ref{tameext}. 
\end{lm}
\proof{
First assume that $\chi$ is unramified. Then $\chi$ comes from an element of ${\rm H}^1(k(\eta)^{\rm ur}/k(\eta),M)$, which is denoted 
by the same letter. Then the character 
$\chi_{[z]}$ coincides with the image of the pair $(z,\chi)\in k(\eta)^\times\times{\rm H}^1(k(\eta)^{\rm ur}/k(\eta),M)$ by the composition of 
\begin{equation*}
k(\eta)^\times\times{\rm H}^1(k(\eta)^{\rm ur}/k(\eta),M)\to{\rm H}^1(k(\eta)^{\rm ur}/k(\eta),k(\eta)^{{\rm ur}\times}\otimes M)
\xrightarrow{{\rm ord}\otimes{\rm id}}{\rm H}^1(k(\eta)^{\rm ur}/k(\eta),M), 
\end{equation*}
where the first arrow is the cup product. Since ${\rm ord}z=1$, the assertion follows. 

Next we treat the general case. For an integer $m\geq1$ which is prime to $p$, denote by $\eta_m$ 
the totally tamely  ramified extension ${\rm Spec}(k(\eta)[X]/(X^m+z))$ of $\eta$. Note that the norm map $k(\eta_m)^\times\to
k(\eta)$ sends $X$ to $z$. 
We use the same letters $\eta$ and $\eta_m$ for their \'etale topoi by abuse of notation. 
Let $\tilde{\eta}$ and $\tilde{\eta_m}$ denote the topoi associated with the \'etale sites 
of unramified extensions of $\eta$  and $\eta_m$ respectively. 
The cohomological functors ${\rm H}^i(\eta_m,-),{\rm H}^i(\tilde{\eta_m},-)$ are naturally identified with 
${\rm H}^i(k(\eta_m)^{\rm sep}/k(\eta_m),-),{\rm H}^i(k(\eta_m)^{\rm ur}/k(\eta_m),-)$ respectively. 
We have a natural commutative diagram 
\begin{equation*}
\xymatrix{
\eta_m\ar[r]^f\ar[d]_{\pi_m}&\eta\ar[d]^\pi\\
\tilde{\eta_m}\ar[r]_{\tilde{f}}&\tilde{\eta}
}
\end{equation*}
of topoi. The norm map $N\colon f_\ast\mathcal{O}_{\eta_m}^\times\to\mathcal{O}_{\eta}^\times$ gives 
the following commutative diagrams 
\begin{equation}\label{comm1}
\xymatrix{
k(\eta_m)^\times\times{\rm H}^1(\eta,M)\ar[r]\ar[d]_{N\times{\rm id}}&k(\eta_m)^\times\times{\rm H}^2(\eta,\mathbb{Z}^{\oplus n})
\ar[d]_{N\times{\rm id}}\ar[r]^{{\rm id}\times f^\ast}&k(\eta_m)^\times\times{\rm H}^2(\eta_m,\mathbb{Z}^{\oplus n})
\ar[r]^{\ \ \ \ \ \ \cup}&
{\rm H}^2(\eta_m,\mathcal{O}^\times_{\eta_m})^{\oplus n}\ar[d]^N\\
k(\eta)^\times\times{\rm H}^1(\eta,M)\ar[r]&k(\eta)^\times\times{\rm H}^2(\eta,\mathbb{Z}^{\oplus n})\ar[rr]^{\cup}&&
{\rm H}^2(\eta,\mathcal{O}^\times_{\eta})^{\oplus n}. 
}
\end{equation}
The commutative diagram 
\begin{equation*}
\xymatrix{
\tilde{f}_\ast\mathcal{O}^\times_{\tilde{\eta_m}}\ar[r]^N\ar[d]&\mathcal{O}^\times_{\tilde{\eta}}\ar[d]\\
\pi_\ast f_\ast\mathcal{O}^\times_{\eta_m}\ar[r]^{\pi_\ast N}&\pi_\ast\mathcal{O}^\times_{\eta}
}
\end{equation*}
induces the commutative diagram 
\begin{equation}\label{comm3}
\xymatrix{
{\rm H}^2(\eta_m,\mathcal{O}^\times_{\eta_m})\ar[d]^N&{\rm H}^2(\tilde{\eta_m},\mathcal{O}^\times_{\tilde{\eta_m}})\ar[l]\ar[d]^N\\
{\rm H}^2(\eta,\mathcal{O}^\times_{\eta})&{\rm H}^2(\tilde{\eta},\mathcal{O}^\times_{\tilde{\eta}}).\ar[l] 
}
\end{equation}
Since $\eta_m/\eta$ is totally ramified, we have a commutative diagram 
\begin{equation}\label{comm4}
\xymatrix{
{\rm H}^2(\tilde{\eta_m},\mathcal{O}^\times_{\tilde{\eta_m}})\ar[r]^{{\rm ord}_{\eta_m}}\ar[d]^N&{\rm H}^2(\tilde{\eta_m},\mathbb{Z})\\
{\rm H}^2(\eta,\mathcal{O}^\times_{\eta})\ar[r]^{{\rm ord}_\eta}&{\rm H}^2(\tilde{\eta},\mathbb{Z})\ar[u]_{\tilde{f}^\ast}
}
\end{equation}
Since the pull-back $\tilde{f}^\ast\colon{\rm H}^1(\tilde{\eta},M)\to{\rm H}^1(\tilde{\eta_m},M)$ is canonically 
identified with the identity map of ${\rm  Hom}(G_k,M)$, the assertion follows by combining (\ref{comm1}), (\ref{comm3}), and (\ref{comm4}). 
\qed}

\begin{lm}\label{gentw}
Let the notation be as above. Let $K\in D_{\rm ctf}(\eta,\Lambda)$ 
be a complex. 
\begin{enumerate}
\item We have 
\begin{equation*}
\varepsilon_{0,\Lambda}(T,K,\omega)\cdot\varepsilon_{0,\Lambda}(T,K,\omega')^{-1}
={\rm det}(K)_{[
\frac{\omega}{\omega'}]}\chi_\cyc^{({\rm ord}(\omega')-{\rm ord}(\omega)){\rm rk}K}. 
\end{equation*}
\item Let $L$ be a locally constant sheaf of finite free $\Lambda$-modules on $T$. 
We have 
\begin{equation*}
\varepsilon_{0,\Lambda}(T,K\otimes^L L,\omega)=
{\rm det}(L)^{a(T,K,\omega)}\cdot
\varepsilon_{0,\Lambda}(T,K,\omega)^{{\rm rk}L}. 
\end{equation*}
The definition of $a(T,K,\omega)$ is given in 
Proposition \ref{elemep}. 
\end{enumerate}
\end{lm}
\proof{
1. It is given in \cite[4.12.(5)]{Y3}. 

2. This follows from the isomorphism 
$F^{(0,\infty)}(K\otimes^L L|_\eta)\cong
L|_\eta\otimes F^{(0,\infty)}(K)$, Theorem \ref{fou}, and $1$. 
\qed
}

\section{Main Results}\label{mr}

In this section, we state and prove the main results. 

First we give notations. Let $T$ be a trait with perfect residue field and 
$f\colon X\to T$ be a morphism of schemes of finite type. Let $\Lambda$ be a 
finite local ring whose residue characteristic is invertible in $T$. 
For a constructible complex $K\in D_{\rm ctf}(X,\Lambda)$ and a closed point 
$x\in  X$ over the closed point $s$ of $T$, we say that $x$ is {\it an at most isolated singularity relative to }$K$ 
if there exists an open neighborhood $U$ of $x$ in $X$ such that 
the restriction $f|_{U\setminus\{x\}}$ is universally locally acyclic relatively to the restriction of $K$. 
Let $T_{(s)}$ and $T_{(x)}$  be the henselization of $T$ and its unramified extension with residue field $k(x)$. 
We denote by 
$\eta$ and $\eta_x$ their generic points respectively. 

Suppose that $x$ is an at most isolated singularity relative to $K$ and $U$ is an open neighborhood of $x$ as above. 
The restriction of 
the vanishing cycles complex $R\Phi_f(K)$ to 
$U\overset\gets\times_TT$ is supported on $x\overset\gets\times_{T_{(s)}}\eta
\cong \eta_x$. 
\begin{df}\label{vancyc}
Let the notation be as above. 
We denote by $R\Phi_f(K)_x$ the pull-back of $R\Phi_f( K)$ by 
$\eta_x\overset\cong\to x\overset\gets\times_{T_{(s)}}\eta\to
U\overset\gets\times_TT$. This is an object of 
$D_{\rm ctf}(\eta_x,\Lambda)$. 
\end{df}
For an object $M\in D_{\rm ctf}(\eta_u,\Lambda)$, define the total dimension 
$\dt M$ of $M$ to be $\dt M:=\rk M+\sw M$.

Let $\Lambda$ be a finite local ring of residual characteristic $\ell\neq p$. 
Fix a non-trivial character $\psi\colon \F_p\to\Lambda^\times$. 
Let $S$ be  a noetherian scheme over $\mathbb{F}_p$. Let 
\begin{equation}\label{rel}
\xymatrix{
Z\ar@{^{(}-_>}[r]&X\ar[rd]_g\ar[rr]^{f}&&Y\ar[ld]\ar[r]^t&\mathbb{A}^1_S\ar[lld] \\
&&S 
}
\end{equation}
be a commutative diagram of $S$-schemes of finite type. 
Let $K\in D_{\rm ctf}(X,\Lambda)$ be a constructible complex of $\Lambda$-sheaves on $X$. 
Consider the following conditions on these data. 
\begin{enumerate}
\item $Y$ is a smooth relative curve over $S$. The morphism $t\colon Y\to\mathbb{A}_S^1$ 
is \'etale. 
\item $Z$ is a closed subscheme of $X$ finite over $S$. 
\item $g$ is universally locally acyclic relatively to $K$. 
\item $f$ is universally locally acyclic relatively to $K$ on 
$X\setminus Z$. 
\end{enumerate}

\begin{lm}\label{bcrel}
Let $X\xrightarrow{f}Y\xrightarrow{g}Z$ be morphisms of finite type of noetherian schemes. 
Let $K\in D_{\rm ctf}(X,\Lambda)$ and $L\in D_{\rm ctf}(Y,\Lambda)$ be constructible complexes on $X$ and $Y$. 
If $f$ (resp. $g$) is universally locally acyclic relatively to $K$ (resp. $L$), so is the composition 
$gf$ relatively to $K\otimes^Lf^\ast L$. 
\end{lm}
\proof{
Let $h_1\colon Z'\to Z$ be a  morphism of noetherian schemes. Let 
\begin{equation*}
\xymatrix{
X'\ar[r]^\beta\ar[d]^{h_3}&Y'\ar[r]^\alpha\ar[d]^{h_2}&Z'\ar[d]^{h_1}\\
X\ar[r]^f&Y\ar[r]^g&Z
}
\end{equation*}
be the base change of $(f,g)$ by $h_1$. By \cite[Theorem 7.6.9]{Fu}, it is enough to show that, for any $M\in D^+(Z',\Lambda)$, the base change map 
\begin{equation}\label{longbc}
(gf)^\ast Rh_{1\ast}M\otimes^L(K\otimes^L f^\ast L)\to Rh_{3\ast}((\alpha\beta)^\ast M\otimes^Lh_{3}^\ast(K\otimes^Lf^\ast L)) 
\end{equation}
defined from the outer square is an isomorphism. 
Note that, since the morphisms $f$ and $g$ are assumed to be of finite type, the universal local acyclicity is 
equivalent to the strongly universal local acyclicity. 
The map (\ref{longbc}) decomposes as follows. 
\begin{align*}
({\rm LHS})\cong
f^\ast(g^\ast Rh_{1\ast}M\otimes^LL)\otimes^LK&\to f^\ast(Rh_{2\ast}(\alpha^\ast M\otimes^Lh_2^\ast L))\otimes^LK\\&
\to Rh_{3\ast}(\beta^\ast(\alpha^\ast M\otimes^Lh_2^\ast L)\otimes^Lh_3^\ast K)\cong({\rm RHS})
\end{align*}
Then the assertion follows from \cite[Theorem 7.6.9]{Fu}. 
\qed}

\begin{lm}\label{rellem}
Let the notation be as above. Assume that the conditions from $1$ to $4$ hold. Then 
the projection $p_2\colon X\times_S\Proj^1_S\to \Proj^1_S$ is universally 
locally acyclic relatively to $p_1^\ast K\otimes^L
((t\circ f)\times{\rm id})^\ast\bar{\mathcal{L}}_\psi$ on $(X\times_S
\Proj^1_S)\setminus(Z\times_S\infty_S)$ 
where $p_1\colon X\times_S\Proj^1_S\to X$ is the projection and $\infty_S
\subset\Proj^1_S$ is the closed subscheme defined by the section 
$S\to\Proj^1_S$ at the infinity.  
\end{lm}
\proof{
Since $t$ is \'etale, we may replace $Y$ and $t$ by $\mathbb{A}^1_S$ and 
the identity map. 
Denote $\mathcal{I}:=p^\ast_1K
\otimes^L(f\times{\rm id})^\ast\bar{\mathcal{L}}_{\psi}$. 
Since $\bar{\mathcal{L}}_{\psi}$ is locally constant on $\mathbb{A}^1_S\times_S
\mathbb{A}^1_S$ and $g$ is universally locally acyclic relatively to $K$, 
$p_2\colon X\times_S\Proj^1_S\to\Proj^1_S$ is universally locally acyclic relatively to $\mathcal{I}$ on $X\times_S
\mathbb{A}^1_S$. Let $U=X\setminus Z$ be the complement. We need to show that 
$U\times_S\Proj^1_S\to\Proj^1_S$ is universally locally acyclic relatively to $\mathcal{I}$. This follows from 
Lemma \ref{bcrel} for the data $U\times_S\Proj^1_S\xrightarrow{f\times{\rm id}}\mathbb{A}^1_S\times_S\Proj^1_S\to\Proj^1_S$, 
$p_1^\ast K\in D_{\rm ctf}(U\times_S\Proj^1_S,\Lambda)$, and $\bar{\mathcal{L}}_\psi\in D_{\rm ctf}(\mathbb{A}^1_S
\times_S\Proj^1_S,\Lambda)$. The universal local acyclicity of the projection $\mathbb{A}^1_{\F_p}\times_{\F_p}\Proj^1_{\F_p}
\to\Proj^1_{\F_p}$ relative to $\bar{\mathcal{L}}_\psi$ is proved in \cite[Th\'eor\`eme (1.3.1.2)]{Lau}. 
\qed}

Before stating Theorem \ref{flat}, we prepare some constructions and 
notation. 

Suppose that the conditions from $1$ to $ 4$ on (\ref{rel}) hold. 
We replace $Y$ and $t$ by $\mathbb{A}^1_S$ and the identity. 
Consider the following diagram of topoi 
\begin{equation*}
\xymatrix{
X\times_S\mathbb{P}^1_S\ar[r]^{\Psi_{f\times{\rm id}}\ \ \ \ \ \ \ \ \ \ \ \ \ \ }\ar[d]_{p_1}
&(X\times_S\mathbb{P}^1_S)
\overset\gets\times_{\mathbb{A}^1_S\times_S\mathbb{P}^1_S}
(\mathbb{A}^1_S\times_S\mathbb{P}^1_S)\ar[r]^{\ \ \ \ \ \ \ \ \ \ \ \ \ q}\ar[d]_{p'}
&\mathbb{A}^1_S\times_S\mathbb{P}^1_S\ar[d]
\\X\ar[r]^{\Psi_f}\ar[rd]_{\rm id}
&X\overset\gets\times_{\mathbb{A}^1_S}\mathbb{A}^1_S
\ar[r]\ar[d]_{{\rm pr}_X}
&\mathbb{A}^1_S\ar[d]^{{\rm id}}
\\&X\ar[r]^f&\mathbb{A}^1_S.
}
\end{equation*}
Consider the following distinguished triangle 
\begin{equation*}
\xymatrix{
{\rm pr}_X^\ast K\ar[r]&
R\Psi_f(K)\ar[r]&
R\Phi_f(K)\ar[r]&
}
\end{equation*}
of complexes on $X\overset\gets\times_{\mathbb{A}^1_S}\mathbb{A}_S$. 
By Proposition \ref{pr}.1, the complexes $R\Psi_f(K)$ and $R\Phi_f(K)$ 
are objects of $D_{{\rm ctf}}(X\overset\gets\times_{\mathbb{A}^1_S}\mathbb{A}^1_S,\Lambda)$ 
and their formations commute with arbitrary base change $S'\to S$. Hence we have a distinguished 
triangle 
\begin{equation*}
\xymatrix{
p'^\ast {\rm pr}_X^\ast K\ar[r]&R\Psi_{f\times{\rm id}}(p_1^\ast K)\ar[r]
&p'^\ast R\Phi_f(K)\ar[r]&
}
\end{equation*}
on $(X\times_S\mathbb{P}^1_S)
\overset\gets\times_{\mathbb{A}^1_S\times_S\mathbb{P}^1_S}
(\mathbb{A}^1_S\times_S\mathbb{P}^1_S)$. 

Let the pull-back of $\bar{\mathcal{L}}_{\psi}$ by $\mathbb{A}_S^1\times_S\mathbb{P}^1_S\to
\mathbb{A}_{\F_p}^1\times_{\F_p}\mathbb{P}^1_{\F_p}$ be denoted by 
the same letter $\bar{\mathcal{L}}_{\psi}$. 
Tensoring $q^\ast\bar{\mathcal{L}}_{\psi}$, we get the 
following distinguished triangle
\begin{equation}\label{tr}
p'^\ast {\rm pr}_X^\ast K\otimes^L
 q^\ast\bar{\mathcal{L}}_{\psi}
\to R\Psi_{f\times{\rm id}}(p_1^\ast K)\otimes^L
 q^\ast\bar{\mathcal{L}}_{\psi}
\to p'^\ast R\Phi_f(K)\otimes^L
 q^\ast\bar{\mathcal{L}}_{\psi}\to. 
\end{equation}
Let $\overset\gets h\colon(X\times_S\mathbb{P}^1_S)
\overset\gets\times_{\mathbb{A}^1_S\times_S\mathbb{P}^1_S}
(\mathbb{A}^1_S\times_S\mathbb{P}^1_S)
\to (X\times_S\mathbb{P}^1_S)\overset\gets\times_{\mathbb{P}^1_S}\mathbb{P}^1_S
$ be the morphism defined from the commutative diagram 
\begin{equation*}
\xymatrix{
X\times_S\mathbb{P}^1_S\ar[r]^{f\times{\rm id}}\ar[rd]_{p_2}&\mathbb{A}^1_S\times_S\mathbb{P}^1_S\ar[d]^h\\
& \mathbb{P}^1_S.
}
\end{equation*}
Taking the push-forward of the triangle (\ref{tr}) through $\overset\gets h$, we get a distinguished triangle 
\begin{multline}\label{5.3}
 \ \overset\gets h_\ast
(p'^\ast {\rm pr}_X^\ast K\otimes^L
 q^\ast\bar{\mathcal{L}}_{\psi})
\to\overset\gets h_\ast(R\Psi_{f\times{\rm id}}(p_1^\ast K)\otimes^Lq^\ast\bar{\mathcal{L}}_{\psi})\\
\to\overset\gets h_\ast(p'^\ast R\Phi_f(K)\otimes^L
 q^\ast\bar{\mathcal{L}}_{\psi})\to
\end{multline}
on $(X\times_S\mathbb{P}^1_S)\overset\gets\times_{\mathbb{P}^1_S}\mathbb{P}^1_S$. 
Here we denote by $\overset\gets h_\ast$ the derived push-forward by abuse of notation.  
Let 
\begin{equation*}
\mathcal{K}\to\mathcal{H}\to\mathcal{G}\to
\end{equation*}
be the distinguished triangle which is the restriction of (\ref{5.3}) to the subtopos 
$(Z\times_S\infty_S)\overset\gets\times_{\mathbb{P}^1_S}\mathbb{A}^1_S\subset(X\times_S\mathbb{P}^1_S)\overset\gets\times_{\mathbb{P}^1_S}\mathbb{P}^1_S$. 
Let $\overset\gets g\colon 
(Z\times_S\infty_S)\overset\gets\times_{\mathbb{P}^1_S}
\mathbb{A}^1_S\to\infty_S
\overset\gets\times_{\mathbb{P}^1_S}\mathbb{A}^1_S$ 
be the morphism of topoi defined from $g|_Z\colon Z\to S$. 

We show that 
the complex $\overset\gets g_\ast\mathcal{G}_n$ is locally constant 
constructible and of finite tor-dimension and that its formation commutes 
with arbitrary base change $S'\to S$. 
To this end, 
we prove the following two claims. 
\begin{lm}\label{acyloc}
\begin{enumerate}
\item The complex $\mathcal{K}$ is acyclic. 
 \item The complex $\overset\gets g_\ast\mathcal{H}$ 
 is locally constant constructible and of finite tor-dimension. 
Its formation commutes with arbitrary base change $S'\to S$. 
 \end{enumerate}
\end{lm}
\proof{
$1$. Take a point $(z\leftarrow t)$ of the topos 
$(Z\times_S\infty_S)\overset\gets\times_{\mathbb{P}^1_S}\mathbb{A}^1_S$. 
Namely, $z$ is a geometric point of $Z$, $t$ is a geometric point of 
$\mathbb{A}^1_S$, and a specialization $\infty_{g(z)}
\leftarrow t$ is given. Here $\infty_{g(z)}$ is the geometric point of $\mathbb{P}^1_S$ defined by 
the compositions 
\begin{equation*}
\xymatrix{
z\ar[r]&Z\ar[r]^g&S\ar[r]^\infty&\mathbb{P}^1_S.
}
\end{equation*}
By (\ref{stalkpush}), the stalk $\mathcal{K}_{(z\leftarrow t)}$ is isomorphic to the following 
 complex: 
 \begin{align*}
 R\Gamma((\mathbb{A}^1_S\times_S\mathbb{P}^1_S)_{(f(z),\infty_{g(z)})}
 \times_{\mathbb{P}^1_{S,(\infty_{g(z)})}}&\mathbb{P}^1_{S,(t)}, 
K_z\otimes^L\bar{\mathcal{L}}_{\psi})\\&\cong R\Gamma((\mathbb{A}^1_S\times_S\mathbb{P}^1_S)_{(f(z),\infty_{g(z)})}
 \times_{\mathbb{P}^1_{S,(\infty_{g(z)})}}\mathbb{P}^1_{S,(t)}, 
 \bar{\mathcal{L}}_{\psi})\otimes^LK_z,
 \end{align*}
 where $K_z$ is a constant complex. Since the projection 
$\mathbb{A}^1_S\times_S\mathbb{P}^1_S\to\mathbb{P}^1_S$ is universally 
locally acyclic relatively to $\bar{\mathcal{L}}_{\psi}$ \cite[Th\'eor\`eme (1.3.1.2)]{Lau} and the restriction of $\bar{\mathcal{L}}_{\psi}$ to 
$\mathbb{A}^1_S\times_S\infty_S$ is zero, 
 we have the assertion. 

$2$. We consider the following cartesian diagrams 
\begin{equation*}
\xymatrix{
X\times_S\mathbb{A}^1_S\ar[d]\ar[r]^{R\Psi'\ \ \ \ \ \ \ \ \ \ \ \ \ \ }&(X\times_S\Proj^1_S)\overset\gets\times_{
\mathbb{A}^1_S\times_S\Proj^1_S}(\mathbb{A}^1_S\times_S\mathbb{A}^1_S)
\ar[d]\ar[r]^{\ \ \ \ \ \ \ \overset\gets h_\ast}&(X\times_S\Proj^1_S)\overset\gets\times_{\Proj^1_S}
\mathbb{A}^1_S\ar[d]\\
X\times_S\mathbb{P}^1_S\ar[r]^{R\Psi_{f\times{\rm id}}\ \ \ \ \ \ \ \ \ \ \ \ \ \ }&(X\times_S\Proj^1_S)\overset\gets\times_{
\mathbb{A}^1_S\times_S\Proj^1_S}(\mathbb{A}^1_S\times_S\mathbb{P}^1_S)
\ar[r]^{\ \ \ \ \ \ \ \overset\gets h_\ast}&(X\times_S\Proj^1_S)\overset\gets\times_{\Proj^1_S}
\Proj^1_S. 
}
\end{equation*}
Since $\bar{\mathcal{L}}_{\psi}$ is locally constant 
on $\mathbb{A}^1_S\times_S\mathbb{A}^1_S$, we compute 
\begin{align}\label{hn}
\mathcal{H}&\cong\overset\gets h_\ast(R\Psi'
(p^\ast_1K)\otimes^Lq^\ast\bar{\mathcal{L}}
_{\psi})|_{(Z\times_S\infty_S)\overset\gets\times_{\mathbb{P}^1_S}\mathbb{A}^1_S}\\&\notag\cong\overset\gets h_\ast(R\Psi'
(p^\ast_1K\otimes^L(f\times{\rm id})^\ast\bar{\mathcal{L}}
_{\psi}))|_{(Z\times_S\infty_S)\overset\gets\times_{\mathbb{P}^1_S}\mathbb{A}^1_S}
\\&\notag\cong R\Psi_{p_2}(p^\ast_1K
\otimes^L(f\times{\rm id})^\ast\bar{\mathcal{L}}_{\psi})|_{(Z\times_S\infty_S)\overset\gets\times_{\mathbb{P}^1_S}\mathbb{A}^1_S}. 
\end{align}
Hence the assertion follows from the conditions $2$, Lemma \ref{rellem}, and Proposition \ref{ula}. 
\qed
}

Let $\mathcal{C}:=\overset\gets g_\ast\mathcal{G}$. 
This is a constructible complex of locally constant $\Lambda$-sheaves
 of finite tor-dimension on 
$\infty_S\overset\gets\times_{\mathbb{P}^1_S}
\mathbb{A}^1_S$. 

\begin{thm}\label{flat}
Suppose that the conditions from $1$ to $ 4$ hold. Then, 
the complex $\mathcal{C}$ on 
$\infty_S\overset\gets\times_{\Proj^1_S}\mathbb{A}^1_S$ constructed above admits the following properties. 
\begin{enumerate}
\item The formation of $\mathcal{C}$ commutes with 
arbitrary base change $S'\to S$. 
\item
Assume that $S={\rm Spec}(k)$ is the spectrum of a perfect field $k$. 
Let $\eta_{k}$ be 
the generic point of the henselization $\Proj^1_{k,(\infty)}$ of $\Proj^1_{k}$ 
at the infinity. 
Then, by the identification $\eta_k\cong\infty\overset\gets\times_{\Proj^1_k}
\mathbb{A}^1_k$, 
$\mathcal{C}$ is isomorphic to 
\begin{equation*}
\bigoplus_{z\in Z}{\rm Ind}_{G_{\eta_z}}^{G_{\eta_{k}}}(
F^{(0,\infty)}(R\Phi_{f}(K)_z)\otimes\mathcal{L}_\psi
(f(z)\cdot x'))[-1],
\end{equation*}
where $\eta_z$ is the unramified extension of $\eta_{k}$ whose residue 
field is isomorphic to that of $z$, $f(z)\in k(z)$ is the image of the standard 
coordinate by $z\to\mathbb{A}^1_{k}$, and $x'$ is the standard coordinate 
of $\mathbb{A}^1_{k(z)}=\Proj^1_{k(z)}\setminus\infty_z$. 
The definition of $R\Phi_{f}(K)_z
$ is given in Definition \ref{vancyc}. 
\end{enumerate}
\end{thm}
\proof{ 

1. This follows from Lemma \ref{acyloc}. 

2. 
Define  
 $g_{(\infty)}\colon\coprod_{z\in Z}\eta_z\to\eta_k$ to be the disjoint union 
of the canonical maps. We have
\begin{align*}
(\overset\gets g_\ast \mathcal{G})_{\overline{\eta_k}}\cong 
 g_{(\infty)\ast}(\mathcal{G}|_{\coprod_z\eta_z})_{\overline{\eta_k}}\cong
\bigoplus_{z\in Z}{\rm Ind}^{G_{\eta_k}}_{G_{\eta_z}}
(\mathcal{G}_{\overline{\eta_z}}). 
\end{align*}
We have 
\begin{align*}
\mathcal{G}_{\overline{\eta_z}}&\cong
R\Gamma((\mathbb{A}^1_{k(z)}\times_{k(z)}\mathbb{P}^1_{k(z)})
_{({f(z)},{\infty_z})}\times_{\mathbb{P}^1_{k(z),({\infty_z})}}
\overline{\eta_z},p'^\ast R\Phi_f(K)\otimes^Lq^\ast
\bar{\mathcal{L}}_{\psi}).
\end{align*}
Thus we get 
\begin{equation*}
\mathcal{C}
\cong \bigoplus_{z\in Z}{\rm Ind}_{G_{\eta_z}}^{G_{\eta_{k}}}(
F^{(0,\infty)}(R\Phi_{f}(K)_z)\otimes\mathcal{L}_\psi
(f(z)\cdot x'))[-1],
\end{equation*}
hence the assertion. 
\qed
}

The inverse of the determinant $\mathop{{\rm det}}\mathcal{C}$ is a locally  constant sheaf of rank $1$ on 
$\infty_S\overset\gets\times_{\Proj^1_S}\mathbb{A}^1_S$, which is 
the tensor product of the determinants of local Fourier transforms and 
Artin-Schreier sheaves. We show that the part of Artin-Schreier sheaves itself 
forms a locally constant sheaf $\mathcal{L}$ of rank $1$. 
\begin{lm}\label{flat-dimtot}
Let the notation be as above. 
Assume that the conditions from $1$ to $4$ are satisfied. 
Let $h\colon Z\to S$ be the structure morphism. 
Then, the map $\varphi_{K}\colon Z\to\mathbb{Z}$ 
defined by $z\mapsto{\rm dimtot}R\Phi_{f_s}(K)_z$, 
where $s$ is the image $h(z)$ of $z$, 
is flat in the sense of Definition \ref{flmoi}. 
\end{lm}
\proof{
This is a consequence of the existence of a locally constant complex 
$\mathcal{C}$ as in Theorem \ref{flat} whose rank 
at $s\in S$ equals to $\sum_{z\in Z_{\bar{s}}}\varphi_{K}
(z)$. It is also proved in \cite[Proposition 2.16]{Sai17}. 
\qed}

Let $\varphi_{K}\colon Z\to \mathbb{Z}$ be the flat function 
defined in Lemma \ref{flat-dimtot}. 
Applying Lemma \ref{AStr} to $\varphi_{K}$ and 
the section $f|_Z\colon Z\to\mathbb{A}^1_Z$ defined from the 
composition $Z\to X\xrightarrow{f}\mathbb{A}^1_S$, 
we obtain a locally constant $\Lambda$-sheaf $\mathcal{L}_\psi(\varphi_{K}\cdot f|_Z)$ of rank $1$ 
on $\mathbb{A}^1_S$. 

\begin{df}
We define the locally constant $\Lambda$-sheaf 
$\mathcal{L}_{\psi,K,f}$ of rank $1$ on $\mathbb{A}^1_S$ to be 
$\mathcal{L}_\psi(\varphi_{K}\cdot f|_Z)$. 
\end{df}

\begin{thm}\label{flatness}
Let the notation be as in Theorem \ref{flat}. Assume that the conditions from $1$ to 
$4$ hold. 
\begin{enumerate}
\item Let ${\rm pr}\colon \infty_S\overset\gets\times_{\Proj^1_S}\mathbb{A}^1_S
\to \mathbb{A}^1_S$ be the second projection. 
Then, the product $\mathop{{\rm det}}(\mathcal{C})^{-1}\otimes {\rm pr}^\ast
\mathcal{L}^{-1}_{\psi,K,f}\in{\rm H}^1(\infty_S
\overset\gets\times_{\Proj^1_S}\mathbb{A}^1_S,\Lambda^\times)$  
is a tame object in the sense of Definition \ref{tameobj1}. 
\item 
Further assume that 
$S$ is connected and normal. Then, 
the continuous group homomorphism
\begin{equation*}
\rho_t\colon\pi_1(S)^{ab}\to \Lambda^\times
\end{equation*}
corresponding to $\langle{\rm det}(\mathcal{C})^{-1}\otimes {\rm pr}^\ast\mathcal{L}_{
\psi,K,f}^{-1},1/x'\rangle\in{\rm H}^1(S,\Lambda^\times)$ defined in Lemma \ref{tameobj} has 
the following properties. 
\begin{enumerate}
\item The formation of $\rho_t$ commutes with arbitrary base change 
$S'\to S$. 
\item When $S={\rm Spec}(k)$ is the spectrum of a perfect field $k$, 
$\rho_t$ 
coincides with 
\begin{equation*}
\prod_{z\in Z}\delta_{k(z)/k}^{{\rm dimtot}(R\Phi_{f}(K)_z)}\cdot\varepsilon_{0,\Lambda}(Y_{(z)},R\Phi_{f}(K)_z,
dt)\circ {\rm tr}_{k(z)/k}. 
\end{equation*}
\end{enumerate}
Here the definition of $\delta_{k(z)/k}$ is given at the end of 
Section \ref{intro}. 
When $S$ is of finite type over $\F_p$, 
$\rho_t$ is uniquely determined by the properties above. 
\item 
Further assume that 
$S={\rm Spec}(k)$ is the spectrum of a perfect field $k$.  
Let $t'\in\Gamma(Y,\mathcal{O}_Y)$ be another section which satisfies 
the conditions $1$ and $4.$ 
We have 
\begin{equation*}
\rho_{t'}\rho_t^{-1}=\prod_{z\in Z}
{\rm det}(R\Phi_{f}(\mathcal{F})_z)_{[\frac{dt'}{dt}]}\circ {\rm tr}_{k(z)/k}. 
\end{equation*}

\end{enumerate}
\end{thm}

\proof{

Let ${\rm Spec}(k)\to S$ be a morphism from the spectrum of a perfect field $k$. 
By Theorem \ref{flat}.1 and 2, we have 
\begin{align}\label{det}
{\rm det}(\mathcal{C})|_{G_k^{ab}}
&\cong\bigotimes_{z\in Z_k}
{\rm det}({\rm Ind}_{G_{\eta_z}}^{G_{\eta_k}}(F^{(0,\infty)}(R\Phi_{f_k}(K_{k})_z
)\otimes\mathcal{L}_\psi(f_k(z)\cdot x')))^{-1}\\
&\notag\cong \mathcal{L}_\psi (-\alpha\cdot x')\otimes\bigotimes_{z\in Z_k}
{\rm det}({\rm Ind}_{G_{\eta_z}}^{G_{\eta_k}}(F^{(0,\infty)}(R\Phi_{f_k}(K_{k})_z
)))^{-1}, 
\end{align}
where $\alpha$ is an element of $k$ defined by $\sum_{z\in Z_k}
\dt R\Phi_{f_k}(K_{k})_z\cdot{\rm Tr}_{k(z)/k}(f_k(z)). $
The subscripts $(-)_k$ mean the base changes by ${\rm Spec}(k)\to S$. We have the equality (\cite[Proposition 1.2]{Del}) 
\begin{equation}\label{eqeqfou}
{\rm det}({\rm Ind}_{G_{\eta_z}}^{G_{\eta_k}}(F^{(0,\infty)}(R\Phi_{f_k}(K_{k})_z
)))=\delta_{k(z)/k}^{{\rm dimtot}(R\Phi_{f_k}(
K_{k})_z)}\cdot{\rm det}(F^{(0,\infty)}(R\Phi_{f_k}(K_{k})_z))\circ {\rm tr}_{k(z)/k}. 
\end{equation}

1.
This follows from (\ref{det}), (\ref{eqeqfou}), and Corollary \ref{tame}. 

2. 
The assertion ($a$) follows from Theorem \ref{flat}.1. 
The assertion ($b$) follows from (\ref{det}) and (\ref{eqeqfou}).  
The last assertion follows from the Chebotarev density. 

3. This follows from Lemma \ref{gentw}.1 and (b).  
\qed
}
\vspace{0.1in}

(Proof of Theorem \ref{eplocep})
Shrinking $Y$ and $S$, we may assume that there is an \'etale morphism 
$t\colon Y\to\mathbb{A}^1_S$ such that $\omega=dt$. The assertion is a special case of Theorem \ref{flatness}.2. 
\qed

\section*{Acknowledgement}
The author would like to express sincere gratitude to his advisor Professor Takeshi Saito for many useful advices and his warm and patient encouragement. This work was supported by the Program for Leading Graduate Schools, MEXT, Japan and by 
JSPS KAKENHI Grant Number 19J11213.


\begin{thebibliography}{99}
\bibitem{Bei}
Beilinson,\ A.:
{\rm Constructible sheaves are holonomic,}
Selecta Math. (N.S.) 22 (2016), no. 4, 1797–1819.
\bibitem{SGA4-3}
Deligne,\ P.:
{\rm Cohomologie \`a supports propres,}
SGA 4 Expos\'e XVII, Th\'eorie des Topos et Cohomologie \'Etale des Sch\'emas, 
Lecture Notes in Mathematics Volume 305, 1973.
\bibitem{SGA7-2}
Deligne,\ P.:
{\rm Le formalisme des cycles \'evanescents,}
SGA 7 II Expos\'e XIII, Groupes de Monodromie en G\'eom\'etrie Alg\'ebrique, 
Lecture Notes in Mathematics Volume 340, 1973.
\bibitem{Del}
Deligne,\ P.:
{\rm Les constantes des équations fonctionnelles des fonctions L,}
Modular functions of one variable, II (Proc. Internat. Summer School, Univ. Antwerp, Antwerp, 1972), pp. 501–597.
\bibitem{Weil2}
Deligne,\ P.:
{\rm La conjecture de Weil. II,}
Inst. Hautes Études Sci. Publ. Math. No. 52 (1980), 137–252. 
\bibitem{Fu}
Fu,\ L.:
{\rm Etale cohomology theory,}
Revised edition. Nankai Tracts in Mathematics, 14. World Scientific Publishing Co. Pte. Ltd., Hackensack, NJ, 2015.


\bibitem{Aby}
Grothendieck,\ A.: 
{\rm Propret\'e cohomologique des faisceaux d'ensembles et des 
faisceaux de groupes non commutatifs,}
SGA 1 Expos\'e XIII, Revetements \'Etales et Groupe Fondamental,
Lecture Notes in Mathematics Volume 224, 1971.
\bibitem{SGA4}
Grothendieck,\ A., Verdier,\ J.- L.:
{\rm Topos,} 
SGA 4 Expos\'e IV, Th\'eorie des Topos et Cohomologie \'Etale des Sch\'emas, Lecture Notes in Mathematics Volume 269, 1972. 
\bibitem{point}
Grothendieck,\ A., Verdier,\ J.- L.:
{\rm Conditions de finitude. Topos et sites fibr\'es. Applications aux questions de passage \`a la limite,}
SGA 4 Expos\'e VI, Th\'eorie des Topos et Cohomologie \'Etale des Sch\'emas, Lecture Notes in Mathematics Volume 270, 1972. 
\bibitem{app}
llusie,\ L.:
{\rm Appendice \`a Th\'eor\`emes de finitude en cohomologie $\ell$-adique},
Cohomologie \'etale SGA 4$\frac12$,
Springer Lecture Notes in Math.\ 569 (1977). 
\bibitem{Ill}
Illusie,\ L.:
{\rm Around the Thom–Sebastiani theorem, with an appendix by Weizhe Zheng,}
Manuscripta Math. 152 (2017), no. 1-2, 61–125.
\bibitem{Autor}
Illusie,\ L.:
{\rm Autour du th\'eor\`eme de monodromie locale, }
P\'eriodes p-adiques (Bures-sur-Yvette, 1988). 
Ast\'erisque No. 223 (1994), 9–57. 
\bibitem{ori}
Illusie,\ L.:
{Exposé XI. Produits orientés,}
Travaux de Gabber sur l'uniformisation locale et la cohomologie étale des schémas quasi-excellents. 
Astérisque No. 363-364 (2014), 213–234.
\bibitem{GK}
Katz,\ N.-M.:
{\rm Local-to-global extensions of representations of fundamental groups,}
Ann. Inst. Fourier (Grenoble) 36 (1986), no. 4, 69–106. 
\bibitem{semi-c}
Laumon,\ G.:
{\rm Semi-continuit\'e du conducteur de Swan (d'apr\`es P. Deligne),}
S\'eminaire E.N.S. (1978-1979) Expos\'e 9, Ast\'erisque 82-83, 173-219 (1981). 
\bibitem {Lau}
Laumon,\ G.:
{\rm Transformation de Fourier, constantes d'équations fonctionnelles et conjecture de Weil,}
Inst. Hautes Études Sci. Publ. Math. No. 65 (1987), 131–210. 
\bibitem{LZ}
Lu,\ Q.,\ Zheng,\ W.:
{\rm Duality and nearby cycles over general bases,}
arXiv:1712.10216. 
\bibitem{Org}
Orgogozo,\ F.:
{\rm Modifications et cycles proches sur une base générale,}
Int. Math. Res. Not. 2006, 1-38.
\bibitem{Sai17}
Saito,\ T.:
{\rm The characteristic cycle and the singular support of a constructible sheaf,}
Inventiones Mathematicae, 207(2) (2017), 597-695.
\bibitem{ep_char}
Takeuchi,\ D.:
{\rm Characteristic epsilon cycles of $\ell$-adic sheaves on varieties}, 
preprint.
\bibitem{Y3}
Yasuda,\ S.:
{\rm Local $\varepsilon_0$-characters in torsion rings,}
J. Théor. Nombres Bordeaux 19 (2007), no. 3, 763–797.
\bibitem{Y1}
Yasuda,\ S.:
{\rm Local constants in torsion rings,}
 J. Math. Sci. Univ. Tokyo 16 (2009), no. 2, 125–197. 
 \bibitem{Y2}
 Yasuda,\ S.:
 {\rm The product formula for local constants in torsion rings,}
 J. Math. Sci. Univ. Tokyo 16 (2009), no. 2, 199–230. 
\end{thebibliography}
\end{document}